%% file: main_imexds.tex
\documentclass[a4paper,11pt]{scrartcl}

\usepackage[USenglish]{babel}
\usepackage[T1]{fontenc}
\usepackage[utf8]{inputenc}

\usepackage{graphicx,caption,subcaption,tikz,pgfplots}
\usepackage[customcolors]{hf-tikz}
\hfsetbordercolor{white}
\definecolor{myblue}{RGB}{17,64,111}
\definecolor{myred}{RGB}{150,30,15}

\usepackage{amsmath,amssymb,amsfonts,amsthm,array,stmaryrd,booktabs,paralist,todonotes,mathtools,siunitx}
\usepackage{array,tablefootnote,thmtools}
\usepackage[toc,page]{appendix}

\usepackage{algorithm,ifthen,placeins}
\usepackage{algpseudocode}

\usepackage{ifpdf}
\ifpdf
\usepackage[pdftex,colorlinks=true,linkcolor=blue,citecolor=green,urlcolor=blue,bookmarks]{hyperref}
\else
\usepackage[dvipdfm,colorlinks=true,linkcolor=blue,citecolor=green,urlcolor=blue,bookmarks]{hyperref}
\fi 

\numberwithin{equation}{section}

\newcommand{\dd}{\textup{d}}
\newcommand{\bfU}{\mathbf{U}}
\newcommand{\bfW}{\mathbf{W}}
\newcommand{\dbl}{\left\llbracket}
\newcommand{\dbr}{\right\rrbracket}

\newcommand{\Sandra}[1]{{\color{red!55!blue}{#1}}}

\begin{document}

    %Titelseite
    \title{Explicit-Implicit Domain Splitting for Two Phase Flows with Phase Transition}
    \author{Sandra May\footnotemark[1],\;Ferdinand Thein\footnotemark[2]}
    \date{\today}
    \maketitle
    \begin{abstract}
        Two phase flows that include phase transition, especially phase creation, with a sharp interface remain a challenging task for numerics.
    We consider the isothermal Euler equations with phase transition between a liquid and a vapor phase.
    The phase interface is modeled as a sharp interface and the mass transfer across the phase boundary is modeled by a kinetic relation.
    Existence and uniqueness results were proven in Ref. \cite{Hantke2019a}.
    %We present a method to obtain the numerical solution for associated Riemann problems.
    Using sharp interfaces for simulating nucleation and cavitation results in the grid containing tiny cells that are several orders of magnitude smaller than the remaining grid cells. This forces explicit time stepping schemes to take tiny time steps on these cells. As a remedy we suggest an explicit implicit domain splitting where the majority of the grid cells is treated explicitly and only the neighborhood of the tiny cells is treated implicitly. We use dual time stepping to solve the resulting small implicit systems. Our numerical results indicate that the new scheme is robust and provides significant speed-up compared to a fully explicit treatment.
    %In particular we show how the cases of nucleation and cavitation may be treated.
    %We will highlight the major difficulties and propose possible strategies to overcome these problems.
        %
    \end{abstract}
    \renewcommand{\thefootnote}{\fnsymbol{footnote}}
    \footnotetext[1]{Department of Information Technology, Uppsala University\\ Box 337, 751 05 Uppsala, Sweden.
    \href{mailto:sandra.may@it.uu.se}{\textit{sandra.may@it.uu.se}}}
    \footnotetext[2]{IGPM, RWTH Aachen,\\ Templergraben 55, D-52056 Aachen, Germany.
    \href{mailto:thein@igpm.rwth-aachen.de}{\textit{thein@igpm.rwth-aachen.de}}}
    \renewcommand{\thefootnote}{\arabic{footnote}}
    \newpage

    %%Inhaltsverzeichnis
    %\tableofcontents
    %\newpage

    %%ToDoListe
    %\listoftodos
    %\newpage
    %%

    %Kapitel
\input{introduction}
%\newpage
\input{num_probl.tex}
%\newpage
\input{num_interface.tex}
%\newpage
%\input{num_nucl_cav.tex}
%\newpage
\input{examples.tex}
%\newpage
\input{conclusion.tex}

%\begin{acknowledgments}
\section*{Acknowledgements}
This work was initiated during a joint participation at the \textit{Hirschegg workshop of conservation laws}. The authors like to thank Marsha Berger for the introduction to dual time stepping.
F.T. is funded by the DFG SPP 2183 \emph{Eigenschaftsgeregelte Umformprozesse}, project 424334423 and gratefully acknowledges the support by the research training group \textit{Energy, Entropy and Dissipative Dynamics (EDDy)} of the DFG - project no. 320021702/GRK2326.
%\end{acknowledgments}

\section*{Data Availability Statement}

No experimental data was produced. Data related to the numerical results is available upon reasonable request.

\appendix
\input{appendix.tex}

%\nocite{*}
\bibliographystyle{abbrv}
\bibliography{literatur_imexds}% Produces the bibliography via BibTeX.
\end{document}

%% file: introduction.tex
\section{Introduction}%\label{sec:intro_2p_iso_num}
In this contribution, we study isothermal liquid-vapor flow problems. In particular, we suggest a new method to compute phase creation phenomena such as cavitation or nucleation.

Due to its diverse applications, the topic of two phase flows is widely discussed in the literature.
Very recently, experimental studies on cavitation have been published that study the underlying dynamics, see Pfeiffer et al.\cite{Pfeiffer2022,Pfeiffer2022a}.
In applications these results can, for example, be used to verify the quality of ultrapure water.
Further usage of cavitation arise in engineering applications, such as rotating under water turbine blades, or in medicine, for example in the treatment of tumor cells.

Mathematically, different models are available to model two phase flows, each with its own benefits and difficulties.
In diffuse interface models such as the Navier-Stokes-Korteweg model, the phase interface is smeared over a certain distance and the interface width is incorporated in the underlying energy potential,
see, e.g., Refs. \cite{Neusser2015,Hitz2020} and the references therein.
So called phase field models have an individual evolution equation for the phase field coupled to the equations governing the flow dynamics, see, e.g., Refs. \ \cite{Dreyer2014,Hantke2022}.
A widely used class models the flow as a diffuse mixture, see Refs. \cite{Baer1986,Saurel1999a,Kapila2001,Romenski2016}. For further reading and a comparative study we also refer to Zein et al.\cite{Zein2010a}.

In the present work we study compressible two phase flows with phase transitions across a sharp interface, see Refs. \cite{Bedeaux2004,Merkle2007,Hantke2013,Fechter2017}.
In particular we consider isothermal inviscid two phase flows governed by the Euler equations, where the phase transition is modeled using a kinetic relation.
This setting was first discussed by Hantke et al.\cite{Hantke2013} and then extended by Thein\cite{Thein2018} and Hantke and Thein\cite{Hantke2019a}.
The use of a kinetic relation for undercompressive shock waves was suggested by Abeyaratne and Knowles \cite{Abeyaratne1991} for solid-solid phase transitions.
This kinetic relation controls the mass transfer across the interface between the two adjacent phases.
For a more general context of kinetic relations see LeFloch\cite{LeFloch2002}.

A fundamental key for a better understanding of hyperbolic problems is the Riemann problem since it exhibits all non-linear phenomena and is widely used as the core ingredient in numerical methods.
A key aspect of Riemann problems is that the constructed solutions are self similar. They consist of constant states, separated by classical rarefaction and shock waves or in the present case of phase boundaries.
Here the considered substance, e.g.,\ water, enters the isothermal Euler equations as the equation of state (EOS) given as a non-monotone pressure-density function.
The pressure function is composed of three parts: the equations of state for the two single phases and an arbitrary relation for the intermediate state.
The two phases are distinguished using the Maxwell construction, also known as the Equal-Area-Rule. The mass transfer is modeled via a kinetic relation, derived in Dreyer et al.\cite{Dreyer2012},
based on classical Hertz-Knudsen theory, see Bond and Struchtrup\cite{Bond2004}.
Without aiming at completeness we refer for further reading to Refs. \cite{Mueller2006,Voss2005} who also considered the isothermal Euler system.
In contrast to the above mentioned work, they model the fluid using the van der Waals equation of state. Instead of a kinetic relation the Liu entropy condition is used to achieve uniqueness.
As a consequence non-classical composite waves are needed to construct solutions.
A further reference for vapor-liquid phase transitions in the context of conservation laws is the detailed review by Fan and Slemrod\cite{Fan2002}.
They also treat the isothermal case using the method of vanishing viscosity applied to the Lagrange formulation of the conservation laws.
Additional literature in this context can be found in the references given before.

An important question is the proper numerical treatment of two phase problems, which of course depends on the chosen mathematical model. The numerical methods in the literature are as diverse as the available models.
Since our focus is on sharp interface methods we briefly refer to Refs. \cite{Dumbser2013,Ancellin2018,Rajkotwala2019,Boniou2021,Bures2021} and the references therein for different methods used in other contexts.
Sharp interface methods are discussed, for example, in Refs. \cite{Schleper2016,Fechter2017,Long2021,Bempedelis2022,Joens2022,Magiera2022}.

Sharp interface methods clearly separate the different phases, and are therefore in better consistency with the underlying physics. 
Also, they avoid the issue of having to find physically meaningful values for the cells in the phase transition layer. For example, when water is considered the density values of the different phases differ by several orders of magnitude, leading to bad results if the averaging in the numerical method
is not performed properly. 

On the technical or numerical side though sharp interface methods pose several challenges. Among other things, one needs to create meshes and keep track of cells, as well as develop methods that are stable on the potentially distorted cells. The higher the dimension, the more complicated this becomes. A particularly challenging situation in the context of two phase flows is the situation of phase creation. 
This is therefore barely discussed in the literature. To make progress in understanding the difficulties and providing suitable numerical methods it is therefore important to first focus on the one dimensional case, which we will do in this contribution. Our goal is to provide a new solution approach for dealing with phase creation in one dimension, which we then plan to extend to higher dimensions in the future.
This method is able to efficiently treat and resolve very small cells that occur during cavitation and nucleation while respecting the sharp interface, and hence avoiding wrong mixing of different phases.

Our approach is based on the mixed explicit implicit time stepping scheme suggested by May and Berger\cite{May_Berger_explimpl}. We here extend it from linear advection to two phase flow problems. The method was originally suggested in the context of Cartesian embedded boundary meshes. This corresponds to a specific way of mesh generation that is suitable to deal with complex geometries. In that approach, a given geometry is cut out of a Cartesian background mesh. Where the object intersects the background mesh, so called \textit{cut cells} are created. These cells can have various shapes and in particular can become arbitrarily small. As a result, when solving time-dependent hyperbolic flow problems on cut cell meshes, one faces the so called \textit{small cell problem} -- that standard explicit time stepping schemes are not stable on the small cut cells if the time step is chosen according to the size of background cells. When dealing with cavitation and nucleation, we experience exactly the same problem. 

Over the years, a number of approaches have been suggested to treat the small cell problem in the context of cut cell meshes for single phase problems, see, e.g., 
Refs. \cite{Chern_Colella,Berger_Helzel_Leveque_2005,Klein_cutcell,FVCA_Helzel_Kerkmann,DoD_SIAM_2020,Kreiss_Fu,Berger_Giuliani_2021,DoD_AMC_2021,Giuliani_DG}. Note that several of the suggested approaches rely on unifying the solution values of the small cell and their bigger neighbors in a suitable way. This is not a feasible option for two phase flow problems. We therefore follow the mixed explicit implicit approach here, where standard explicit time stepping is used away from the small cell and implicit time stepping is used in the neighborhood of the small cell for stability. A similar approach of combining explicit and implicit time stepping has been used very recently by Fu et al.\cite{Fu_interface} in the context of locally moving interfaces for the advection problem for higher-order discontinuous Galerkin methods. 

As a result of the mixed time stepping, one has to solve an implicit system involving the small cells and its neighbors in each time step. In the context of phase transitions this involves complicated non-linear Riemann solvers. Differentiating through them, as needed for a Newton approach, is non-trivial. Therefore, in this work we use so called dual time stepping to solve the implicit systems. In that approach, which we will describe in more detail in section \ref{sec:num_interface} below, the work that needs to be done in each iteration of the implicit solver is very similar to the work of taking an explicit time step.

This contribution is organized in the following way.
First we will introduce the problem setting and point out the main challenges in section \ref{sec:num_probl}.
In section \ref{sec:num_interface}, we will introduce the mixed explicit implicit scheme as well as the dual time stepping routine for solving the resulting implicit systems.
In section \ref{sec:num_ex} we will present numerical results for our new method including tests for cavitation and nucleation.
We conclude with a short summary in section \ref{sec:conclusion}.
%

%% file: num_probl.tex
\section{Problem Formulation and Main Challenges}\label{sec:num_probl}
\subsection{Description of two phase flows}
For the description of the flow dynamics under study in the present work we will follow Hantke and Thein\cite{Hantke2019a} and give a brief summary in the following.\\
We study inviscid, compressible and isothermal two phase flows.
The two phases are either the liquid or the vapor phase of one substance. The phases are distinguished by the \textit{mass density} $\rho$, and are further described by the \textit{velocity} $u$.
Sometimes it is convenient to use the \textit{specific volume} $v = 1/\rho$ instead of the mass density.
%We will make the reader aware of such situations.
The physical quantities depend on time $t \in \mathbb{R}_{\geq 0}$ and space $x \in \mathbb{R}$.

In regular points of the bulk phases the fluid is described using the (one dimensional) isothermal Euler equations
\begin{subequations}\label{eeq:isoT both}
\begin{align}
    \partial_t\rho + \partial_x(\rho u) &= 0,\label{eeq:isoT:cons_mass_1d}\\
    \partial_t(\rho u) + \partial_x(\rho u^2 + p) &= 0.\label{eeq:isoT:cons_mom_1d}
\end{align}
\end{subequations}
The role of the mathematical entropy inequality is here played by the energy inequality
\begin{align}
    0 &\geq -T\zeta = \frac{\partial}{\partial t}\left(\rho f + \rho\frac{u^2}{2}\right)
    + \partial_x\left(\rho u\left(g + \frac{u^2}{2}\right)\right).\label{eeq:iso_T:entropy_balance}
\end{align}
Here $f$ is the specific free energy and $g = f + p/\rho$ is the Gibbs energy of the bulk phase.
This is the proper stability condition for a thermodynamic system with uniform and constant temperature, see Refs. \cite{Mueller1985,Landau1991}.
Note that this inequality also includes the heat flux, which is needed as a mechanism for the isothermal process.

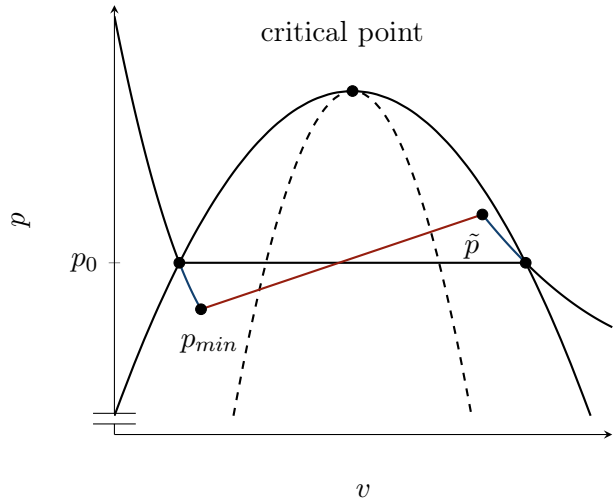
\begin{figure}[!ht]
    \center
    \begin{tikzpicture}
        \begin{axis}[y=1cm,x=1cm,xlabel={$v$},axis x line=bottom,ylabel={$p$},axis y line=left,axis y discontinuity = parallel, ymin=0,ymax=5,domain=0:10,axis equal
        image,ytick={2},yticklabel={$p_0$},xtick=\empty]
            %saturation curve
            \addplot[mark=none,smooth,black,thick,solid,domain=0.25:5.75] (\x,{-0.5*(\x-3)^2 + 4});
            %spinodal curve
            \addplot[smooth,black,thick,dashed,domain=1.625:4.375] (\x,{-2*(\x-3)^2 + 4});
            %saturation pressure
            \addplot[smooth,black,thick,solid,domain=1:5] (\x,2);
            %mark points
            \addplot[mark=*] coordinates {(1, 2)};
            \addplot[mark=*] coordinates {(3, 4)};
            \addplot[mark=*] coordinates {(5, 2)};
            \addplot[mark=*] coordinates {(1.25, 1.459040820577346)};
            \addplot[mark=*] coordinates {(4.5, 41/16)};
            %isotherm
            \addplot[mark=none,smooth,black,thick,solid,domain=0.25:1] (\x,{1.663836717690617*\x^2 - 5.907469332494504*\x + 6.243632614803888});
            \addplot[mark=none,smooth,myblue,thick,solid,domain=1:1.25] (\x,{1.663836717690617*\x^2 - 5.907469332494504*\x + 6.243632614803888});
            \addplot[mark=none,smooth,myred,thick,solid,domain=1.25:4.5] (\x,{0.339525901360817*(\x - 1.25) + 1.459040820577346});
            \addplot[smooth,myblue,thick,solid,domain=4.5:5] (\x,{0.25*(\x-7)^2+1});
            \addplot[smooth,black,thick,solid,domain=5:6] (\x,{0.25*(\x-7)^2+1});

        \end{axis}
        %
        %mark critical point
        \draw[color=black] (3,5) node[above] {critical point};
        \draw[color=black] (1.25, 1.459040820577346) node[below] {$p_{min}$};
        \draw[color=black] (4.7, 2.8) node[below] {$\tilde{p}$};
    \end{tikzpicture}
    \caption{Schematic figure of the $v-p$ phase plane with the EOS satisfying the Maxwell construction.}
    \label{fig:maxwell_construct}
\end{figure}

The pressure $p$ is linked to the mass density $\rho$ via the EOS. For the mass density we have $\rho \in \Omega_\rho \subseteq (0,\infty)$.
This domain can be split into the vapor, spinodal, and liquid region, i.e.,\ $\Omega_\rho = \Omega_{vap}\cup\Omega_{spin}\cup\Omega_{liq}$ with
\[
  \Omega_{vap} := (0,\tilde{\rho}],\;\Omega_{spin} := (\tilde{\rho},\rho_{\min}),\; \Omega_{liq} := [\rho_{\min},\infty).
\]
Accordingly the EOS consists of three corresponding parts, i.e.,\ an EOS for the vapor phase, the liquid phase, and an intermediate part, see figure \ref{fig:maxwell_construct}.
The EOS satisfies the Maxwell construction (equal-area rule). For further details we refer to Refs. \cite{Thein2018,Hantke2019a}.
From the Maxwell construction we can obtain the maximum vapor pressure $\tilde{p}$ and the minimum liquid pressure $p_{min}$.
In the regular phases the EOS has the properties discussed in Ref. \cite{Hantke2019a}. We also note that the (unphysical) intermediate part is characterized by the relation
\begin{align*}
    \left(\frac{\partial p}{\partial v}\right)_T > 0.
\end{align*}
Thus the considered Euler system becomes elliptic inside this region.
It is this region that is crossed when phase transitions are present leading to the failure of the standard theory for hyperbolic systems.
Here the state space is separated into two regions where the system is hyperbolic and these may be connected by so called undercompressive shock waves.

In the present work phase boundaries are considered as sharp interfaces and thus are treated as discontinuities, in particular as undercompressive shocks.
These undercompressive shock waves are not of Lax type and hence are lacking the usual uniqueness properties.
However, uniqueness can be restored by prescribing an additional algebraic relation. Although this choice remains non-unique there is a preferred choice based on physical considerations.
For further reading we refer to Refs. \cite{Dafermos2016,LeFloch2002,Merkle2007}.

Since we have a conservative system given by \eqref{eeq:isoT both}, the following jump conditions hold across discontinuities
\begin{subequations}
\begin{align}
    \dbl \rho(u - W)\dbr &= 0,\label{jc:isoT:cons_mass_1d}\\
    \rho(u - W)\dbl u\dbr + \dbl p\dbr &= 0.\label{jc:isoT:cons_mom_1d}
\end{align}
\end{subequations}
Here we write $[\![\Psi]\!] = \Psi^+ - \Psi^-$, where $\Psi^+$ is the right and $\Psi^-$ the left sided limit of the physical quantity $\Psi$.
Furthermore every discontinuity satisfies the following entropy inequality
\begin{align}
    \rho(u - W)\dbl g + e^{kin}\dbr \leq 0,\label{isoT:entropy_ineq_1d}
\end{align}
which is consistent with \eqref{eeq:iso_T:entropy_balance}.
Here $g$ denotes the specific Gibbs energy and $e_{kin}$ the specific kinetic energy.
The quantity $W$ is the speed of the discontinuity and $Z = -\rho(u - W)$ is the mass flux, where we will distinguish between a classical shock wave and the phase boundary (non-classical shock)
\begin{align*}
    Z = \begin{cases} Q,\;\text{shock wave,}\\ z,\;\text{phase boundary,}\end{cases}
    \text{and}\;
    W = \begin{cases} S,\;\text{shock wave,}\\ w,\;\text{phase boundary.}\end{cases}
\end{align*}
Our choice of the kinetic relation, which defines the mass transfer across the interface, is the following, compare Refs. \cite{Dreyer2012,Hantke2019a},
\begin{align}
    z = \tau p_V\dbl g + e^{kin}\dbr.\label{isoT:kin_rel1}
\end{align}
Here $p_V$ denotes the pressure of the vapor phase and $\tau > 0$.
We verify that this kinetic relation satisfies \eqref{isoT:entropy_ineq_1d}
\begin{align*}
    \rho(u - W)\dbl g + e^{kin}\dbr = -z\dbl g + e^{kin}\dbr = -\tau p_V\dbl g + e^{kin}\dbr^2 \leq 0.
\end{align*}

\subsection{Riemann problems}

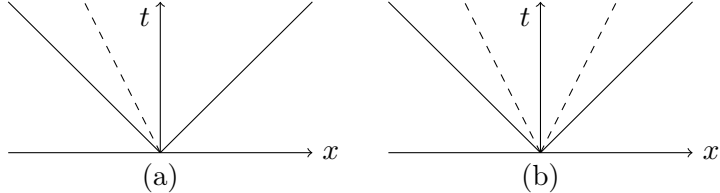
\begin{figure}[h!]%Different wave patterns
    \begin{center}
        \begin{tikzpicture}[scale=1,domain=-7:7]
            %
            %2phase init data
            %axes
            \draw[->] (-7,0) -- (-3,0) node[right] {$x$};
            \draw (-5,0) node[below] {(a)};
            \draw[->] (-5,0) -- (-5,2);
            \draw[color=black] (-5,1.8) node[left] {$t$};
            %left wave
            \draw[color=black] (-5,0) -- (-7,2);
            %middle wave
            \draw[color=black,dashed] (-5,0) -- (-6,2);
            %right wave
            \draw[color=black] (-5,0) -- (-3,2);
            %
            %single phase init data
            %axes
            \draw[->] (-2,0) -- (2,0) node[right] {$x$};
            \draw (0,0) node[below] {(b)};
            \draw[->] (0,0) -- (0,2);
            \draw[color=black] (0,1.8) node[left] {$t$};
            %left wave
            \draw[color=black] (0,0) -- (-2,2);
            %middle wave
            \draw[color=black,dashed] (0,0) -- (-1,2);
            \draw[color=black,dashed] (0,0) -- (1,2);
            %right wave
            \draw[color=black] (0,0) -- (2,2);
        \end{tikzpicture}
    \end{center}
    \caption{Wave patterns. Solid line: classical waves. Dashed line: phase boundary. }
    \label{diff_wave_patterns}
\end{figure}

A particular choice of initial data is the Riemann initial data. For the isothermal Euler equations for two phases with and without phase transition the Riemann problem was discussed in Refs. \cite{Hantke2013,Hantke2019a}.
Thus we will only briefly review the main aspects.
The Riemann initial data is given by
\begin{align}
    \rho(x,0) = \begin{cases} \rho^-,\, &x < 0,\\ \rho^+,\,&x > 0,\end{cases}
    \quad\text{and}\quad
    u(x,0) = \begin{cases} u^-,\, &x < 0,\\ u^+,\,&x > 0,\end{cases}\label{isoT:init_data_rp}
\end{align}
where the constant states left and right of the initial discontinuity may belong to different phases or to the same phase which then may lead to nucleation or cavitation.
The solution of the Riemann problem is self similar and in the case of two phase initial data it consists of three waves separating four constant states where the phase boundary lies in between the classical outer waves, see figure \ref{diff_wave_patterns} (a). We provide more information about how to solve the Riemann problem between two phases in appendix \ref{sec:app}.

For the case of single phase initial data with phase creation the solution consists of four waves separating five constant states where the phase boundaries lie in between the classical outer waves, see figure \ref{diff_wave_patterns} (b).
Across each wave specific relations hold that are used to obtain values inside the wave fan. We provide additional information about phase creation in appendix \ref{sec:num_cav_nucl}.
For further details we again refer to the previous mentioned Refs. \cite{Hantke2013,Hantke2019a}.

\subsection{Sharp interfaces and grid alignment}

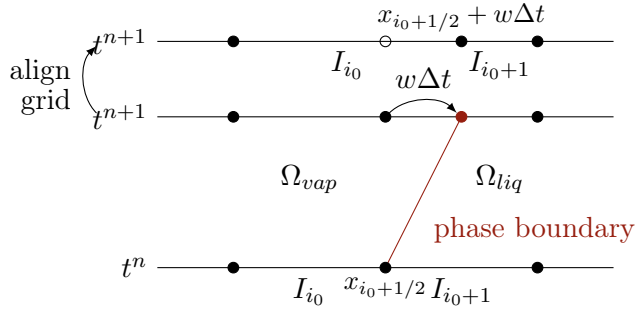
\begin{figure}[h!]
    \hspace*{-1.75cm}
    %\begin{subfigure}[]{0.15\textwidth}
    \centering
    \begin{tikzpicture}
        %
        %phase boundary
        \draw[color=myred] (3,0) -- (4,2);
        \draw[color=myred] (3.5,0.5) node[right] {phase boundary};
        \draw[color=black] (2,1.5) node[below] {$\Omega_{vap}$};
        \draw[color=black] (4.5,1.5) node[below] {$\Omega_{liq}$};
        %space axis at t
        \draw[color=black] (0,0) node[left] {$t^n$};
        \draw[color=black] (0,0) -- (6,0);
        \draw[mark=*] plot coordinates {(1,0)};
        \draw[mark=*] plot coordinates {(3,0)};
        \draw[mark=*] plot coordinates {(5,0)};
        \draw[color=black] (2,0) node[below] {$I_{i_0}$};
        \draw[color=black] (3,0) node[below] {\small$x_{i_0 + 1/2}$\normalsize};
        \draw[color=black] (4,0) node[below] {$I_{i_0+1}$};

        %space axis at t + Delta t
        \draw[color=black] (0,2) node[left] {$t^{n+1}$};
        \draw[color=black] (0,2) -- (6,2);
        \draw[mark=*] plot coordinates {(1,2)};
        \draw[mark=*] plot coordinates {(3,2)};
        \draw[mark=*,color=myred] plot coordinates {(4,2)};
        \draw[mark=*] plot coordinates {(5,2)};
        \draw[-latex,shorten <= 2pt,shorten >= 2pt] (3,2) arc [radius=0.7, start angle = 135, end angle = 45];
        \draw[color=black] (3.5,2.25) node[above] {$w\Delta t$};
        %aligned grid
        \draw[-latex,shorten <= 2pt,shorten >= 2pt] (-0.75,2) arc [radius=0.7, start angle = 225, end angle = 135];
        \draw[color=black] (-1,2.6) node[left] {align};
        \draw[color=black] (-1,2.2) node[left] {grid};
        \draw[color=black] (0,3) node[left] {$t^{n+1}$};
        \draw[color=black] (0,3) -- (6,3);
        \draw[mark=*] plot coordinates {(1,3)};
        \draw[mark=o] plot coordinates {(3,3)};
        \draw[mark=*] plot coordinates {(4,3)};
        \draw[mark=*] plot coordinates {(5,3)};
        \draw[color=black] (2.5,3) node[below] {$I_{i_0}$};
        \draw[color=black] (4,3) node[above] {\small$x_{i_0 + 1/2} + w\Delta t$\normalsize};
        \draw[color=black] (4.5,3) node[below] {$I_{i_0+1}$};
\end{tikzpicture}
    \caption{Align grid to the new position of the phase boundary.}
    \label{fig:align_grid_to_pb}
\end{figure}

When doing numerical simulations based on a finite volume approach for two phase flows using sharp interfaces, one faces the following issue:
%Generally, when using sharp interfaces for treating two phase flows the key difficulty when it comes to numerical simulations is the following: 
%We will apply a finite volume approach to compute numerical approximations.
By definition, the unknowns $\bfU_i^n$ correspond to cell averages in cell $I_i$ at time $t^n$. Let us assume that at time $t^n$ our grid is aligned to the phase boundaries. We focus on a single boundary moving at the velocity $w > 0$, located at $x_{i_0 + 1/2}$ at $t^n$, see also figure \ref{fig:align_grid_to_pb}. If we keep the mesh fixed, then at time $t^{n+1}$ the cell $I_{i_0 + 1}$ would contain both partially vapor and partially liquid. Computing an unknown $\bfU_{i_0+1}^{n+1}$ by averaging would most likely result in values inside $\Omega_{spin}$ as well as in a diffused interface. 

Instead we adjust the mesh appropriately by moving the cell boundary $x_{i_0+1/2}$ to 
$x_{i_0 + 1/2} + w\Delta t$ as indicated in figure \ref{fig:align_grid_to_pb}. As a result, the cell size changes within the time step. 
Every now and then cells become too small or too big. We then merge the cells with neighboring cells or split them in two parts.
%We will provide more details in section \ref{subsec: grid alignment}. 

\begin{figure}%[]{0.15\textwidth}
    \centering
    \begin{tikzpicture}
        %
        %phase boundary
        \draw[color=myred] (3,0) -- (3.15,2);
        \draw[color=myred] (3,0) -- (2.85,2);
        \draw[color=myred] (3.3,1) node[right] {right phase boundary};
        \draw[color=myred] (2.7,1) node[left] {left phase boundary};
        %\draw[color=black] (2,1.5) node[below] {$\Omega_{vap/liq}$};
        %\draw[color=black] (4.5,1.5) node[below] {$\Omega_{vap/liq}$};
        %space axis at t
        \draw[color=black] (0,0) node[left] {$t^n$};
        \draw[color=black] (0,0) -- (6,0);
        \draw[mark=*] plot coordinates {(1,0)};
        \draw[mark=*] plot coordinates {(3,0)};
        \draw[mark=*] plot coordinates {(5,0)};
        \draw[color=black] (2,0) node[below] {$I_{i_0}$};
        \draw[color=black] (3,0) node[below] {\small$x_{i_0 + 1/2}$\normalsize};
        \draw[color=black] (4,0) node[below] {$I_{i_0+1}$};

        %space axis at t + Delta t
        \draw[color=black] (0,2) node[left] {$t^{n+1}$};
        \draw[color=black] (0,2) -- (6,2);
        \draw[mark=*] plot coordinates {(1,2)};
        \draw[mark=*,color=myred] plot coordinates {(2.85,2)};
        \draw[mark=*,color=myred] plot coordinates {(3.15,2)};
        \draw[mark=*] plot coordinates {(5,2)};
        %\draw[-latex,shorten <= 2pt,shorten >= 2pt] (3,2) arc [radius=0.7, start angle = 135, end angle = 45];
        \draw[color=black] (3,2.25) node[above] {$(w_{right} - w_{left})\Delta t$};
\end{tikzpicture}
    \caption{For initial states that belong to one phase phase creation (i.e.,\ cavitation or nucleation) can occur. In the simulation this results in the creation of a small cell from one time step to the next. This cell must not be merged with its neighbors since they belong to different phases.}
    \label{fig:phase_creation}
\end{figure}
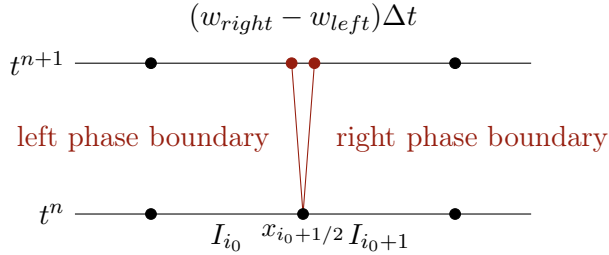

Note that it is \textit{only} possible to merge small cells with one of their neighbors if they belong to the same phase. In the situation where a new phase is created, see figure \ref{fig:phase_creation}, this is \textit{not} possible. We therefore need numerical algorithms that can deal with such meshes that contain tiny cells. 
%

%% file: num_interface.tex
\section{Time stepping in numerical method}\label{sec:num_interface}

In this section we present the explicit implicit time stepping approach. Our goal is to treat situations of phase creation, as shown in figure \ref{fig:phase_creation}, efficiently. 

\subsection{Time stepping and tiny cells}

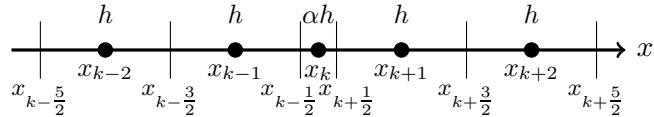
\begin{figure}[htb]
\begin{center}
\begin{tikzpicture}[scale=0.95,
axis/.style={very thick, line join=miter, ->}]
\draw [axis] (-4.0,0) -- (4.5,0) node(xline)[right] {$x$};
\draw (-3.6,-0.4) -- (-3.6,0.4);
\draw (-1.8,-0.4) -- (-1.8,0.4);
\draw (0,-0.4) -- (0,0.4);
\draw (0.5,-0.4) -- (0.5,0.4);
\draw (2.3,-0.4) -- (2.3,0.4);
\draw (4.1,-0.4) -- (4.1,0.4);
  \draw[fill] (-2.7,0) circle (.1cm);
\node[] at (-2.7,-0.32) {\small $x_{k-2}$};
  \draw[fill] (-0.9,0) circle (.1cm);
\node[] at (-0.9,-0.32) {\small $x_{k-1}$};
  \draw[fill] (0.25,0) circle (.1cm);
\node[] at (0.25,-0.32) {\small $x_k$};
  \draw[fill] (1.4,0) circle (.1cm);
\node[] at (1.4,-0.32) {\small $x_{k+1}$};
  \draw[fill] (3.2,0) circle (.1cm);
\node[] at (3.2,-0.32) {\small $x_{k+2}$};
\node[] at (-3.6,-0.7) {\footnotesize $x_{k-\tfrac 5 2}$};
\node[] at (-1.8,-0.7) {\footnotesize $x_{k-\tfrac 3 2}$};
\node[] at (-0.15,-0.7) {\footnotesize $x_{k-\tfrac 1 2}$};
\node[] at (0.65,-0.7) {\footnotesize $x_{k+\tfrac 1 2}$};
\node[] at (2.3,-0.7) {\footnotesize $x_{k+\tfrac 3 2}$};
\node[] at (4.1,-0.7) {\footnotesize $x_{k+\tfrac 5 2}$};
\node[] at (-2.7,0.5) {\small $h$};
\node[] at (-0.9,0.5) {\small $h$};
\node[] at (0.25,0.5) {\small $\alpha h$};
\node[] at (1.4,0.5) {\small $h$};
\node[] at (3.2,0.5) {\small $h$};
\end{tikzpicture}
\end{center}
\caption{1d model problem: equidistant mesh of mesh width $h$ with one small cell (cell $I_k$) of length $\alpha h$, $\alpha \in (0,1]$. Cell $I_k$ corresponds to the creation of a new phase.}
\label{fig:1d model problem} 
\end{figure}

As a simplified model problem for explaining our approach, we consider the mesh shown in figure \ref{fig:1d model problem}. For simplicity, we assume that all other cells have the same length $h$ but this is not necessary. 
% The main point here is that there will be one tiny cell and the other cells have roughly similar cell lengths. 
Also, for now we ignore the fact that the length of the small cell might change slightly during the time step and assume the mesh to be fixed. We denote the cell centroid of cell $I_i$ with $x_i$ and the edges with $x_{i\pm \frac 1 2}$. One can think of the model mesh \ref{fig:1d model problem} to correspond to the mesh at time $t^{n+1}$ in figure \ref{fig:phase_creation} and to keep that mesh fixed.

In space we use a first-order finite volume scheme, also known as \textit{Godunov's method}. 
Our unknowns $\bfU_i^n$ approximate the cell averages of the true solution in cell $I_i$ with length $h_i$ at time $t^n$, i.e., 
\begin{equation}
    \mathbf{U}_i^n \approx \frac{1}{h_i}\int_{x_{i-\frac{1}{2}}}^{x_{i+\frac{1}{2}}}\mathbf{U}(t^n,x)\,\dd x.\label{def:cell_av}
\end{equation}
To update the solution from time $t^n$ to $t^{n+1}$ we need numerical fluxes at cell edges. Denoting the flux at cell boundary $x_{i \pm \frac 1 2}$ with $\mathbf{F}_{i \pm \frac{1}{2}}$, the update is then given by
%\Sandra{Die 2.zeile in der folgenden Formel passt noch nicht in die story}
\begin{equation}\label{godunov_meth:unif}
    \mathbf{U}^{n+1}_i = \mathbf{U}^n_i - \frac{\Delta t}{h_i}\left[\mathbf{F}_{i+\frac{1}{2}} - \mathbf{F}_{i-\frac{1}{2}}\right].
\end{equation}
%\begin{align}
%    %
%    \mathbf{U}^{n+1}_i &= \mathbf{U}^n_i - \frac{\Delta t}{h_i}\left[\mathbf{F}_{i+\frac{1}{2}} - \mathbf{F}_{i-\frac{1}{2}}\right],\label{godunov_meth:unif}.
%    \mathbf{F}_{i+\frac{1}{2}} &= \mathbf{F}\left(\mathbf{U}_{i+\frac{1}{2}}\left(0\right)\right).\label{godunov_flux}
%    %
%\end{align}
Generally, the fluxes $\mathbf{F}_{i + \frac 1 2}$ and $\mathbf{F}_{i - \frac 1 2}$ are computed by solving Riemann problems at cell boundary $x_{i+\frac 1 2}$ with input data $\bfU_i$ and $\bfU_{i+1}$ and at cell boundary $x_{i-\frac 1 2}$ with input data $\bfU_{i-1}$ and $\bfU_{i}$, respectively. 
Here, for fluxes at interfaces between cells with equal phases a standard Riemann solver can be used, see e.g. Toro\cite{Toro2009} for an overview. We specify our choice when presenting the numerical results. For the cell boundaries with phase transition, i.e., the interfaces between cells $I_{k-1}$ and $I_k$ as well as cells $I_k$ and $I_{k+1}$, we employ special Riemann solvers. For completeness, we provide the description in the appendix \ref{sec:app}. The focus of this contribution is on the choice of the time stepping scheme, i.e., on which instance of time to evaluate the input arguments $\bfU_{i \pm 1}$ and $\bfU_{i}$.

The by far most common approach is to employ explicit time stepping. When using explicit Euler in time, this results in
\begin{equation}\label{eq: update expl time}
\mathbf{U}^{n+1}_i = \mathbf{U}^n_i - \frac{\Delta t}{h_i}\left[\mathbf{F}_{i+\frac{1}{2}}(\bfU_{i}^n,\bfU_{i+1}^n) - \mathbf{F}_{i-\frac{1}{2}}(\bfU_{i-1}^n,\bfU_{i}^n)\right].
\end{equation}
Note that the data $\bfU_{i \pm 1}^n$ and $\bfU_{i}^n$ is known and can then be used in the Riemann solver as input. 
For stability, the time step must be chosen to satisfy the \emph{CFL-condition}
\begin{equation}
    \Delta t \leq C_{CFL}\frac{h_i}{S^n_{max}}.\label{cfl_cond:unif}
\end{equation}
Here, $S^n_{max}$ denotes the maximum absolute wave speed in every time step during the numerical simulation. In particular, we expect $\Delta t$ to change in the course of the simulation, i.e., to be precise it should be $\Delta t_n$. For better readability, we will drop that index though.
The constant $C_{CFL}$ is the \emph{CFL-number}. Typically, $C_{CFL} \in (0,1]$.

If we want to take the \textit{same} time step length for \textit{all} cells, we must choose $\Delta t$ according to our smallest cell. 
Defining $\Delta t_{\text{ref}} = C_{CFL}\frac{h}{S^n_{max}}$, we get for the model grid in 
figure \ref{fig:1d model problem} the CFL constraint
$
 \Delta t \leq \alpha \Delta t_{\text{ref}}.
$
For the nucleation problem that we solve below in section \ref{subsubsec:ex3}, $\alpha=\mathcal{O}(10^{-7})$, forcing us to take a factor of $\mathcal{O}(10^{7})$ more time steps than needed without the presence of the small cell $I_k$. Roughly speaking we expect the running time of the algorithm to also grow with this factor.

One option to reduce some of this additional work is to take a big time step of length $\Delta t_{\text{ref}}$ on cells away from the small cell $I_k$ and to do some form of time-accurate local time stepping in the neighborhood of the small cell $I_k$, i.e., on cells $I_{k-1}, I_k, I_{k+1}$. Full details for this approach can be found with M\"uller and Stiriba \cite{Mueller2007}. For our situation, it breaks down to the following: First, one uses a standard explicit Euler time step with step length $\Delta t_{\text{ref}}$ to update 
 all cells with indices $i \le k-2$  and $i\ge k+2$. Then, the values on cells $I_{k-1},I_k,I_{k+1}$ are updated: here, one applies a time-accurate local time stepping using a step length $\Delta \tau = \mathcal{O}(\alpha \Delta t_{\text{ref}})$.
% (assuming that the maximum wave speed does not change in the time interval $(t^n,t^{n+1})$). 
 That means that one needs to take $\mathcal{O}(\frac{1}{\alpha})$ local time steps for each big time step of length $\Delta t_{\text{ref}}.$ 

Another approach is to use implicit time stepping as this does not require such a strict CFL condition for stability. Implicit Euler in particular has very good stability properties. % and is for example unconditionally strong stability preserving. 
Using implicit Euler results in  an update formula of the form
%Therefore, one option is to use implicit Euler time stepping everywhere resulting in an update formula of the form
\begin{equation}\label{eq: update impl time}
\mathbf{U}^{n+1}_i = \mathbf{U}^n_i - \frac{\Delta t}{h_i}\left[\mathbf{F}_{i+\frac{1}{2}}(\bfU_{i}^{n+1},\bfU_{i+1}^{n+1}) - \mathbf{F}_{i-\frac{1}{2}}(\bfU_{i-1}^{n+1},\bfU_{i}^{n+1})\right].
\end{equation}
Note that now the input arguments $\bfU_{i \pm 1}^{n+1}$ and $\bfU_{i}^{n+1}$ are generally not known. Instead, they are exactly the solution we are looking for. Therefore some kind of implicit solver is necessary in each time step. 

For accuracy reasons, we still would want to choose the time step length $\Delta t$ proportional to the (average) cell length. We could for example choose $\Delta t = \Delta t_{\text{ref}}$, i.e., choose $\Delta t$ based on the size of the bigger cells. The big advantage compared to using explicit time stepping is that we would still be stable on the small cell $I_k$ in this case.
The downside is that the method is much more expensive. In each step we now would need to solve a big non-linear implicit system. 

\subsection{Explicit implicit domain splitting}

We therefore suggest a mixed explicit implicit approach: only treat the neighborhood of the small cell $I_k$ implicit for stability but treat all the cells away from the small cell explicit to keep the cost low. The question is how to change between explicit and implicit time stepping while ensuring conservation and maintaining stability.
May and Berger\cite{May_Berger_explimpl} examined these questions in the context of cut cell meshes. There one also faces similar meshes as shown in figure \ref{fig:1d model problem} that contain tiny cells. The authors suggest to use \textit{flux bounding} to couple the explicit and the implicit scheme. The work in Ref. \cite{May_Berger_explimpl} only considered the linear advection equation. Here we extend this to more complicated non-linear flow problems, in particular to two phase flow problems. 

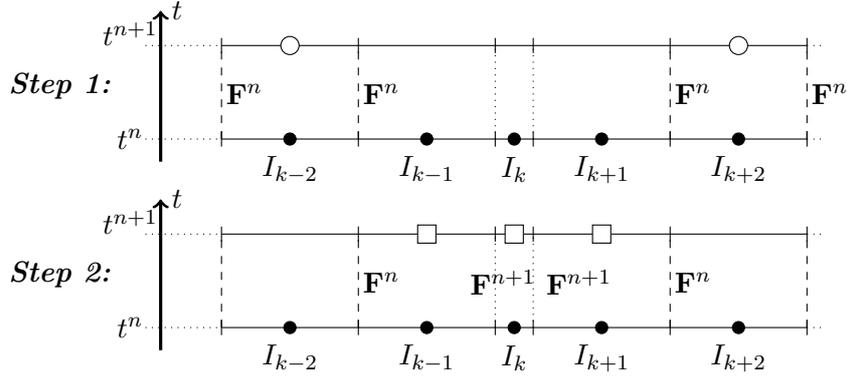
\begin{figure}[htb]
\begin{center}
\begin{tikzpicture}[scale=1,%0.85,
axis/.style={very thick, line join=miter, ->}]
\draw [axis] (-4.4,-0.3) -- (-4.4,1.7) node(xline)[right] {$t$};
\node[left] at (-4.9,0.7) {\textbf{\em Step 2:}};
\node[] at (-4.8,0) {$t^n$};
\node[] at (-4.8,1.4) {$t^{n+1}$};
\draw[dotted] (-4.6,0) -- (-3.6,0);
\draw (-3.6,0) -- (4.1,0);
\draw[dotted] (4.1,0) -- (4.3,0);
\draw (-3.6,-0.1) -- (-3.6,0.1);
\draw (-1.8,-0.1) -- (-1.8,0.1);
\draw (0,-0.1) -- (0,0.1);
\draw (0.5,-0.1) -- (0.5,0.1);
\draw (2.3,-0.1) -- (2.3,0.1);
\draw (4.1,-0.1) -- (4.1,0.1);
\node[] at (-2.7,-0.4) {$I_{k-2}$};
\node[] at (-0.9,-0.4) {$I_{k-1}$};
\node[] at (0.25,-0.4) {$I_k$};
\node[] at (1.4,-0.4) {$I_{k+1}$};
\node[] at (3.2,-0.4) {$I_{k+2}$};
\draw[fill] (-2.7,0) circle (.08cm);
\draw[fill] (-0.9,0) circle (.08cm);
\draw[fill] (0.25,0) circle (.08cm);
\draw[fill] (1.4,0) circle (.08cm);
 \draw[fill] (3.2,0) circle (.08cm);
\draw[dotted] (-4.6,1.24) -- (-3.6,1.24);
\draw (-3.6,1.24) -- (4.1,1.24);
\draw[dotted] (4.1,1.24) -- (4.3,1.24);
\draw (-3.6,1.14) -- (-3.6,1.34);
\draw (-1.8,1.14) -- (-1.8,1.34);
\draw (0,1.14) -- (0,1.34);
\draw (0.5,1.14) -- (0.5,1.34);
\draw (2.3,1.14) -- (2.3,1.34);
\draw (4.1,1.14) -- (4.1,1.34);
%
%\draw[fill=white] (-2.7,1.24) circle (.12cm);
\draw[fill=white] (-1.02,1.12) rectangle (-0.78,1.36);
% \draw (-0.9,1.24) circle (.12cm);
\draw[fill=white] (0.13,1.12) rectangle (0.37,1.36);
%\draw (0.25,1.24) circle (.12cm);
\draw[fill=white] (1.28,1.12) rectangle (1.52,1.36);
%\draw (1.4,1.24) circle (.12cm);
%\draw[fill=white] (3.2,1.24) circle (.12cm);
%
\draw[dashed] (-3.6,-0) -- (-3.6,1.24);
\draw[dashed] (-1.8,-0) -- (-1.8,1.24);
\draw[dotted] (0,-0) -- (0,1.24);
\draw[dotted] (0.5,-0) -- (0.5,1.24);
\draw[dashed] (2.3,-0) -- (2.3,1.24);
\draw[dashed] (4.1,-0) -- (4.1,1.24);
%
%\node[] at (-3.3,0.6){$F^E$};
\node[] at (-1.5,0.6){$\mathbf{F}^n$};
\node[] at (0.1,0.6){$\mathbf{F}^{n+1}$};
\node[] at (1.1,0.6){$\mathbf{F}^{n+1}$};
\node[] at (2.6,0.6){$\mathbf{F}^n$};
%\node[] at (4.4,0.6){$F^E$};
%
%==================================
%
%
\draw [axis] (-4.4,2.2) -- (-4.4,4.2) node(xline)[right] {$t$};
\node[left] at (-4.9,3.2) {\textbf{\em Step 1:}};
\node[] at (-4.8,2.5) {$t^n$};
\node[] at (-4.8,3.9) {$t^{n+1}$};
\draw[dotted] (-4.6,2.5) -- (-3.6,2.5);
\draw (-3.6,2.5) -- (4.1,2.5);
\draw[dotted] (4.1,2.5) -- (4.3,2.5);
\draw (-3.6,2.4) -- (-3.6,2.6);
\draw (-1.8,2.4) -- (-1.8,2.6);
\draw (0,2.4) -- (0,2.6);
\draw (0.5,2.4) -- (0.5,2.6);
\draw (2.3,2.4) -- (2.3,2.6);
\draw (4.1,2.4) -- (4.1,2.6);
\node[] at (-2.7,2.1) {$I_{k-2}$};
\node[] at (-0.9,2.1) {$I_{k-1}$};
\node[] at (0.25,2.1) {$I_k$};
\node[] at (1.4,2.1) {$I_{k+1}$};
\node[] at (3.2,2.1) {$I_{k+2}$};
\draw[fill] (-2.7,2.5) circle (.08cm);
\draw[fill] (-0.9,2.5) circle (.08cm);
\draw[fill] (0.25,2.5) circle (.08cm);
\draw[fill] (1.4,2.5) circle (.08cm);
 \draw[fill] (3.2,2.5) circle (.08cm);
\draw[dotted] (-4.6,3.74) -- (-3.6,3.74);
\draw (-3.6,3.74) -- (4.1,3.74);
\draw[dotted] (4.1,3.74) -- (4.3,3.74);
\draw (-3.6,3.64) -- (-3.6,3.84);
\draw (-1.8,3.64) -- (-1.8,3.84);
\draw (0,3.64) -- (0,3.84);
\draw (0.5,3.64) -- (0.5,3.84);
\draw (2.3,3.64) -- (2.3,3.84);
\draw (4.1,3.64) -- (4.1,3.84);
\draw[fill=white] (-2.7,3.74) circle (.12cm);
\draw[fill=white] (3.2,3.74) circle (.12cm);
\draw[dashed] (-3.6,2.5) -- (-3.6,3.74);
\draw[dashed] (-1.8,2.5) -- (-1.8,3.74);
\draw[dotted] (0,2.5) -- (0,3.74);
\draw[dotted] (0.5,2.5) -- (0.5,3.74);
\draw[dashed] (2.3,2.5) -- (2.3,3.74);
\draw[dashed] (4.1,2.5) -- (4.1,3.74);
\node[] at (-3.3,3.1){$\mathbf{F}^n$};
\node[] at (-1.5,3.1){$\mathbf{F}^n$};
%\node[] at (0.2,3.1){$F^I$};
%\node[] at (0.8,3.1){$F^I$};
\node[] at (2.6,3.1){$\mathbf{F}^n$};
\node[] at (4.4,3.1){$\mathbf{F}^n$};
\end{tikzpicture}
\end{center}
\caption{Idea behind flux bounding for the time step $t^n \to t^{n+1}$ (compare Ref. \cite{FVCA_May}): \textbf{\em Step 1:} All cells away from the small cell $I_k$ are updated (indicated by the symbol `$\circ$') using a standard explicit scheme based on the explicit flux $\mathbf{F}^n$; \textbf{\em Step 2:} The neighborhood of the small cell is updated (indicated by the symbol `$\square$') using implicit fluxes $\mathbf{F}^{n+1}$ for the faces $x_{k \pm \frac 1 2}$.}% Note that in both steps, the same fluxes must be used across edges $x_{k\pm \frac 3 2}$.}
\label{fig:flux bounding} 
\end{figure}

We explain the general idea using again the model mesh
from figure \ref{fig:1d model problem}. To keep the implicit region as small as possible, we first update all cells that are not direct neighbors of cell $I_k$ using explicit time stepping as given by \eqref{eq: update expl time}. This is illustrated in figure \ref{fig:flux bounding}, \textbf{\em Step 1}. For stability, we need a fully implicit update on the small cell 
$I_k$ as given in \eqref{eq: update impl time}. The critical point is how to choose the fluxes $\mathbf{F}_{k-\frac 3 2}$ and $\mathbf{F}_{k+ \frac 3 2}$. When we updated the solution on cells $I_{k-2}$ and $I_{k+2}$ in \textbf{\em Step 1}, we used the fluxes $\mathbf{F}_{k-\frac 3 2}(\bfU_{k-2}^n,\bfU_{k-1}^n)$ and $\mathbf{F}_{k+\frac 3 2}(\bfU_{k+1}^n,\bfU_{k+2}^n)$. We need to reuse these same fluxes for the cell updates for cells $I_{k-1}$ and $I_{k+1}$ in \textbf{\em Step 2} as otherwise conservation may be violated. We therefore get the following updates in the neighborhood of the small cell $I_{k}$
\begin{subequations}\label{eq: update mixed EI}
\begin{align}
 \mathbf{U}^{n+1}_{k-1} &= \mathbf{U}^n_{k-1} - \frac{\Delta t}{h}\left[\mathbf{F}_{k-\frac{1}{2}}(\bfU_{k-1}^{n+1},\bfU_{k}^{n+1}) - \mathbf{F}_{k-\frac{3}{2}}(\bfU_{k-2}^{n},\bfU_{k-1}^{n})\right],\\
    \mathbf{U}^{n+1}_k &= \mathbf{U}^n_k - \frac{\Delta t}{\alpha h}\left[\mathbf{F}_{k+\frac{1}{2}}(\bfU_{k}^{n+1},\bfU_{k+1}^{n+1}) - \mathbf{F}_{k-\frac{1}{2}}(\bfU_{k-1}^{n+1},\bfU_{k}^{n+1})\right],\\
     \mathbf{U}^{n+1}_{k+1} &= \mathbf{U}^n_{k+1} - \frac{\Delta t}{h}\left[\mathbf{F}_{k+\frac{3}{2}}(\bfU_{k+1}^{n},\bfU_{k+2}^{n}) - \mathbf{F}_{k+\frac{1}{2}}(\bfU_{k}^{n+1},\bfU_{k+1}^{n+1})\right].
\end{align}
\end{subequations}
Note how the update on the small cell $I_{k}$ is fully implicit by using fluxes $\mathbf{F}_{k+\frac{1}{2}}(\bfU_{k}^{n+1},\bfU_{k+1}^{n+1})$ and 
$\mathbf{F}_{k-\frac{1}{2}}(\bfU_{k-1}^{n+1},\bfU_{k}^{n+1})$ and how solution updates on cells $I_{k \pm 1}$ have an explicit flux on one side and an implicit flux on the other side. This is also illustrated in figure \ref{fig:flux bounding}, \textbf{\em Step 2}.

To sum it up, for a model mesh like shown in figure \ref{fig:1d model problem}, in each time step,  
\begin{enumerate}
    \item we first update all cells with indices $i \le k-2$ and $i\ge k+2$ using a fully explicit update as given in \eqref{eq: update expl time};
    \item we then update solution values on cells $I_{k-1}, I_k, I_{k+1}$ using \eqref{eq: update mixed EI}.
\end{enumerate}
We note that for a fully explicit treatment using time-accurate local time stepping on the neighborhood of the small cell, the first step would be the same. The difference is how to update the small cell and its face neighbors.

By construction, the mixed explicit implicit scheme preserves conservation. Further, it was shown in Ref.\ \cite{May_Berger_explimpl} that for the linear advection equation $u_t + c u_x = 0$ ($c$ constant) on a model mesh as given in figure \ref{fig:1d model problem}, the scheme is monotonicity preserving and total variation diminishing for a time step length $\Delta t \le \frac{h}{\lvert c \rvert}$, i.e., independent of $\alpha$. In other words, this way of coupling preserves desirable stability properties. In terms of accuracy, numerical tests for linear advection show first order accuracy for the mixed explicit implicit scheme as described here. 

\subsection{Dual time stepping}\label{subsec: dual time stepping}

In the mixed explicit implicit approach, the implicit time stepping 
only couples solution values on cells $I_{k-1},I_k,I_{k+1}$. Nevertheless, given that the update formulae involve non-linear flux functions, in particular Riemann solvers for phase changes, solving the resulting implicit system is non-trivial. 

We suggest dual time stepping for doing this. Dual time stepping can be interpreted as an iterative method, for which the work in each iteration looks very similar to executing a time step with an explicit time stepping scheme. We expect a more Newton-like approach for solving the implicit system to be more efficient and will explore that in the future. The advantage of dual time stepping is that there is no need to differentiate through the Riemann solver. The implementation work for changing a running code from explicit time stepping to implicit time stepping that uses dual time stepping to solve the resulting implicit systems is quite small.

Dual time stepping, or also called pseudo time stepping, has already been around for a number of years
and has become more popular again recently, see, e.g., Refs. \cite{Jameson_AIAA91,Jameson_Sriram,Birken_Jameson,LOPPI2019}. 
The idea is the following: assume we have a first-order finite volume space discretization and use implicit Euler in time for the update on cell $I_i$ resulting in the update formula given in \eqref{eq: update impl time}. We introduce a function $G$ as
\[%begin{multline*}
G(\bfW) = \frac{1}{\Delta t} \left( - \bfW_i + \bfU_i^n  \right) \\- \frac{1}{h_i} \left[\mathbf{F}_{i+\frac{1}{2}}(\bfW_{i},\bfW_{i+1}) - \mathbf{F}_{i-\frac{1}{2}}(\bfW_{i-1},\bfW_{i})\right].
\]%end{multline*}
Obviously, if we find $\bfW$ such that $G(\bfW) = 0$, then $\bfW_i$ corresponds to the searched solution $\bfU_i^{n+1}$. 
In order to solve the implicit system described by $G(\bfW)=0$, we reinterpret this as the steady state solution of the time-dependent problem 
\[
\partial_{\tau} \bfW - G(\bfW) = 0,
\]
i.e., we introduce a dual or pseudo time scale with respect to $\tau$. We can now solve for the steady state by using a standard time stepping scheme for doing time stepping and by running until the solution does not change anymore, i.e., until we have reached (computational) steady state. We use explicit Euler time stepping for the arguments $\bfW$ within the flux evaluation and implicit Euler time stepping for the isolated instance of $\bfW_i$ resulting in the time step formula
\[%begin{multline*}
\bfW_i^{l+1} = \bfW_i^l + \frac{\Delta \tau}{\Delta t} \left( - \bfW_i^{l+1} + \bfU_i^n  \right) \\- \frac{\Delta \tau}{h_i} \left[\mathbf{F}_{i+\frac{1}{2}}(\bfW_{i}^l,\bfW_{i+1}^l) - \mathbf{F}_{i-\frac{1}{2}}(\bfW_{i-1}^l,\bfW_{i}^l)\right].
\]%end{multline*}
Here, $\bfW^{l}$ denotes the given iterate at time $\tau^l$ and $\bfW^{l+1}$ the iterate at $\tau^{l+1} = \tau^l + \Delta \tau$ that we want to compute. We solve for $\bfW^{l+1}$ to get
\begin{align}\label{eq: dual time tau equal}
\bfW_i^{l+1} = \frac{1}{1+\frac{\Delta \tau}{\Delta t}} \left( \bfW_i^l + \frac{\Delta \tau}{\Delta t}   \bfU_i^n
- \frac{\Delta \tau}{h_i} \left[\mathbf{F}_{i+\frac{1}{2}}(\bfW_{i}^l,\bfW_{i+1}^l) - \mathbf{F}_{i-\frac{1}{2}}(\bfW_{i-1}^l,\bfW_{i}^l)\right] \right).
\end{align}
Note that the work for executing this iteration is essentially the same as taking one explicit Euler step for the original scheme, i.e., for evaluating \eqref{eq: update expl time}. In particular, everything on the right hand side is known and can simply be evaluated. On the downside, as we treat the fluxes explicitly when running to steady state, we need to enforce the CFL condition \eqref{cfl_cond:unif}. As we use dual time stepping to compute the updates in the neighborhood of the small cell, this results in the time step constraint
$ 
\Delta \tau \leq \alpha \Delta t_{\text{ref}}.
$

%\begin{algorithm}[htp]
\begin{table*}[htp]
\caption{Algorithms 1 and 2 for implementing the mixed explicit implicit scheme using dual time stepping for mesh in figure \ref{fig:1d model problem}.}\label{tab: algorithms}
\hrule~\newline\vspace*{-8pt}\hrule
%\hline\hline
\vspace*{0.2cm}

\textbf{Algorithm 1:} Executing 1 time step for the explicit implicit scheme.\vspace*{0.1cm}

\begin{algorithmic}
\Require $\bfU^n$, $\Delta t$, $\alpha$, $h$
\State \textbf{\em Step 1:} Update all cells away from small cell using explicit time stepping:
\State \quad (a) Based on $\bfU^n$, compute fluxes $\mathbf{F}$ using a standard Riemann solver.
\State \quad (b) On cells $I_i$ with $i \le k{-}2$ or $i\ge k{+}2$, update using
\[
\bfU_{i}^{n+1} = \bfU_{i}^{n} - \frac{\Delta t}{h} \left( \mathbf{F}_{i+\frac 1 2}(\bfU_{i}^n,\bfU_{i+1}^{n}) - \mathbf{F}_{i-\frac 1 2}(\bfU_{i-1}^{n},\bfU_{i}^n)  \right).
\]
\State \qquad\quad Save fluxes
$  \mathbf{F}_{\text{left}} = \mathbf{F}_{k-\frac 3 2}(\bfU_{k-2}^{n},\bfU_{k-1}^n)$ and
 $  \mathbf{F}_{\text{right}} = \mathbf{F}_{k+\frac 3 2}(\bfU_{k+1}^n,\bfU_{k+2}^{n})$.
\State \textbf{\em Step 2:} Compute the update on cells $I_{k-1}$, $I_k$, $I_{k+1}$ using algorithm 2.
\Ensure $\bfU^{n+1}$ 
\end{algorithmic}
\vspace*{0.2cm}
\hrule
%\hline
\vspace*{0.2cm}
\textbf{Algorithm 2:} Compute update on cells $I_{k-1},I_k,I_{k+1}$ using dual time stepping.\vspace*{0.1cm}

\begin{algorithmic}
\Require $\bfU^n_{k-1}$, $\bfU^n_k$, $\bfU^n_{k+1}$, $\Delta t$, $\alpha$, $h$, $\mathbf{F}_{\text{left}}$, $\mathbf{F}_{\text{right}}$, $TOL$, $l_{\text{max}}$, $\theta_{k-1},\theta_k,\theta_{k+1}$.
\State \textbf{Part 1:} Compute approximations for $\bfU^{n+1}_{k-1}$, $\bfU^{n+1}_k$, $\bfU^{n+1}_{k+1}$:
\State \quad \textit{Notation:} $\bfU^{n,l}$: iterate $l$ in dual time stepping loop for approximating $\bfU^{n+1}$
\State \quad Initialize $\bfU_i^{n,0} = \bfU_i^n$ for $i=k{-}1,k,k{+}1$
\State \quad Set local time steps: $\Delta \tau_i = \theta_i \Delta t$ for $i=k{-}1,k,k{+}1$
%$\Delta \tau_{k-1} = 0.9 \Delta t$, $\Delta \tau_k = \alpha \Delta t$, $\Delta t_{k+1} = 0.9 \Delta t$
\For{$l = 0,\dots,l_{\text{max}}$ }
\State (a) Compute fluxes $\mathbf{F}_{k-\frac 1 2}(\bfU_{k-1}^{n,l},\bfU_{k}^{n,l})$ and $\mathbf{F}_{k+\frac 1 2}(\bfU_{k}^{n,l},\bfU_{k+1}^{n,l})$ using a Riemann solver with phase transition.
\State (b) Update
{\small
\begin{align*}
\bfU_{k-1}^{n,l+1} &= \frac{\bfU_{k-1}^{n,l} + \frac{\Delta \tau_{k-1}}{\Delta t} \left( \bfU_{k-1}^n - \frac{\Delta t}{h} \left( \mathbf{F}_{k-\frac 1 2}(\bfU_{k-1}^{n,l},\bfU_{k}^{n,l}) -   \mathbf{F}_{\text{left}}\right) \right)}{1+\frac{\Delta \tau_{k-1}}{\Delta t}},\\
\bfU_k^{n,l+1} &=\frac{\bfU_{k}^{n,l} +  \frac{\Delta \tau_k}{\Delta t} \left(\bfU_k^n - \frac{\Delta t}{\alpha h} \left( \mathbf{F}_{k+\frac 1 2 }(\bfU_{k}^{n,l},\bfU_{k+1}^{n,l}) - \mathbf{F}_{k-\frac 1 2}(\bfU_{k-1}^{n,l},\bfU_{k}^{n,l})  \right) \right)}{1+\frac{\Delta \tau_k}{\Delta t}},\\
\bfU_{k+1}^{n,l+1} &= \frac{\bfU_{k+1}^{n,l} + \frac{\Delta \tau_{k+1}}{\Delta t} \left(\bfU_{k+1}^n -  \frac{\Delta t}{h} \left( \mathbf{F}_{\text{right}} - \mathbf{F}_{k+\frac 1 2}(\bfU_{k}^{n,l},\bfU_{k+1}^{n,l})  \right) \right)}{1+\frac{\Delta \tau_{k+1}}{\Delta t}}.
\end{align*}
}
\If {$\text{\em Stopping criterium } < TOL$}
\State $\bfU^{n,*}_{i} = \bfU^{n,l+1}_{i}$ for $i=k{-}1,k,k{+}1$
\State Break \quad \% leave for loop 
\EndIf
%\State (d) Update $\bfU_i^{n,l+1} = \bfU_i^{n,l}$ for $i=k-1,k,k+1$, $l = l+1$
\EndFor
\State \textbf{Part 2:} Compute $\bfU^{n+1}_{k-1}$, $\bfU^{n+1}_k$, $\bfU^{n+1}_{k+1}$ based on
$\bfU^{n,*}_{k-1}$, $\bfU^{n,*}_k$, $\bfU^{n,*}_{k+1}$:
\State \quad (a) Compute fluxes $\mathbf{F}_{k-\frac 1 2}(\bfU_{k-1}^{n,*},\bfU_k^{n,*})$ and $\mathbf{F}_{k+\frac 1 2}(\bfU_{k}^{n,*},\bfU_{k+1}^{n,*})$ using a Riemann solver with phase transition.
\State \quad (b) Apply the conservative update formulae
{\small 
\begin{align*}
\bfU_{k-1}^{n+1} &= \bfU_{k-1}^n - \frac{\Delta t}{h} \left( \mathbf{F}_{k-\frac 1 2}(\bfU_{k-1}^{n,*},\bfU_{k}^{n,*}) -   \mathbf{F}_{\text{left}}\right), \\
\bfU_k^{n+1} &= \bfU_k^n - \frac{\Delta t}{\alpha h} \left( \mathbf{F}_{k+\frac 1 2}(\bfU_{k}^{n,*},\bfU_{k+1}^{n,*}) - \mathbf{F}_{k-\frac 1 2}(\bfU_{k-1}^{n,*},\bfU_{k}^{n,*})  \right), \\
\bfU_{k+1}^{n+1} &= \bfU_{k+1}^n -  \frac{\Delta t}{h} \left( \mathbf{F}_{\text{right}} - \mathbf{F}_{k+\frac 1 2}(\bfU_{k}^{n,*},\bfU_{k+1}^{n,*})  \right).
\end{align*}
}
\Ensure $\bfU^{n+1}_{k-1}$, $\bfU^{n+1}_k$, $\bfU^{n+1}_{k+1}$
\end{algorithmic}
\hrule~\newline\vspace*{-8pt}\hrule
%\hline\hline
\end{table*}
%\end{algorithm}

So what did we gain by going this complicated way to dual time stepping? The big advantage is that we do \textit{not need to be time-accurate} anymore when running to the steady state. We are allowed to use local time stepping where we update each cell based on their own local CFL and therefore direct neighbors might advance at different speeds. Therefore, the correct formulation of \eqref{eq: dual time tau equal} is
\begin{align}\label{eq: dual time tau diff}
\bfW_i^{l+1} = \frac{1}{1+\frac{\Delta \tau_i}{\Delta t}} \left( \bfW_i^l + \frac{\Delta \tau_i}{\Delta t}   \bfU_i^n 
- \frac{\Delta \tau_i}{h_i} \left[\mathbf{F}_{i+\frac{1}{2}}(\bfW_{i}^l,\bfW_{i+1}^l) - \mathbf{F}_{i-\frac{1}{2}}(\bfW_{i-1}^l,\bfW_{i}^l)\right] \right),
\end{align}
with using individual $\Delta \tau_i$ instead of a global $\Delta \tau$.

In our model problem, we treat three cells (namely $I_{k-1},I_k,I_{k+1}$) implicitly with cell lengths $h, \alpha h, h$ (and $\alpha \ll 1$). Using the notation $\Delta \tau_i = \theta_i \Delta t_{\text{ref}}$, we need to enforce 
$\theta_{k\pm1} \in (0,1]$ and $\theta_{k} \in (0,\alpha]$ for stability (assuming that the maximum wave speed during the time $t^n \to t^{n+1}$ does not change). We comment on our specific choices of $\theta_{k \pm 1}$ and $\theta_k$ in section
\ref{sec:num_ex}, where we present numerical results.

Note that there is a huge difference to using time-accurate local time stepping as sketched above in a fully explicit setup. First, for time-accurate local time stepping, for the update of the solution on cell $I_k$ from $t^n$ to $t^{n+1}$, one needs to take $\mathcal{O}(\frac{1}{\alpha})$ time steps. % (assuming that the maximum wave speed does not change). 
In our approach here we only take as many iterations as needed to reach steady state. Second, in the time-accurate local time stepping, for the update on the small cell one always needs to extract the correct flux data from the bigger neighbors by a sort of interpolation. This is also not needed here. 

We summarize all the necessary steps for implementing dual time stepping in our setting in algorithms 1 and 2 in table \ref{tab: algorithms}. Note that in algorithm 2, when not executing \textbf{Part 2} but simply using the output $\bfU^{n,*}$ as approximation for $\bfU^{n+1}$, one might violate conservation. Therefore, in \textbf{Part 2} fluxes are recomputed based on the accepted iterate and the standard conservative update formula is applied. We comment on our choices for the stopping criteria and for $\theta_i$ in section \ref{sec:num_ex}.

% \FloatBarrier

\subsection{Accounting for changing cell sizes}\label{subsec: grid alignment}
Most of the discussion so far focused on the model mesh in figure \ref{fig:1d model problem} with fixed cell lengths. As mentioned above, when using sharp interfaces, cell lengths can change. 
Shifting the cell boundaries results in cells becoming smaller or larger. Every now and then we need to reorganize by splitting cells that are too large (e.g., more than twice the length of an average cell) in two parts or by unifying two cells if one gets too small. The latter is of course only possible if the two cells belong to the same phase. For details we refer to Thein \cite{Thein2018}.

Thus we need to take into account that the cell sizes can change during the time step of length $\Delta t$. Generally, we expect the tiny cell that contains the new phase to become bigger and therefore its neighbors to become a bit smaller. 
Therefore, instead of using
\eqref{godunov_meth:unif}, the general update formula needs to be changed to 
\begin{equation}
    \mathbf{U}^{n+1}_i = \frac{h_i^n}{h_i^{n+1}}\mathbf{U}^n_i - \frac{\Delta t}{h_i^{n+1}}\left[\mathbf{F}_{i+\frac{1}{2}} - \mathbf{F}_{i-\frac{1}{2}}\right],
    \label{godunov_meth:nonunif}
\end{equation}
with $h_i^n$ and $h_i^{n+1}$ being the cell sizes at the times $t^n$ and $t^{n+1}$.
One immediately verifies that (\ref{godunov_meth:nonunif}) reduces to (\ref{godunov_meth:unif}) for constant cell sizes. Note that when using explicit time stepping, $h_i^{n+1}$ is predicted explicitly by using wave speed information from the Riemann solutions of 
$\mathbf{F}_{i+\frac{1}{2}}(\bfU_i^n,\bfU_{i+1}^n)$ and $\mathbf{F}_{i-\frac{1}{2}}(\bfU_{i-1}^n,\bfU_{i}^n)$.

All formulae given above for the mixed explicit implicit approach need to be adjusted accordingly. For example, the update for $\bfU_k^{n,l+1}$ in algorithm 2 then is given by
%\[
%\bfU_k^{n,l+1} =\frac{\bfU_{k}^{n,l} +  \frac{\Delta \tau_k}{\Delta t} \left(\frac{h_k^n}{h_k^{n,l+1}}\bfU_k^n - \frac{\Delta t}{h_k^{n,l+1}} %\left( \mathbf{F}_{k+\frac 1 2 }(\bfU_{k}^{n,l},\bfU_{k+1}^{n,l}) - \mathbf{F}_{k-\frac 1 2}(\bfU_{k-1}^{n,l},\bfU_{k}^{n,l})  \right) %\right)}{1+\frac{\Delta \tau_k}{\Delta t}}
%\]
\[%\begin{multline*}
\bfU_k^{n,l+1} = \frac{\bfU_{k}^{n,l}}{1+\frac{\Delta \tau_k}{\Delta t}} +  \frac{\frac{\Delta \tau_k}{\Delta t}}{1+\frac{\Delta \tau_k}{\Delta t}}  \left(\frac{h_k^n}{h_k^{n,l+1}}\bfU_k^n \right. \\
\left. - \frac{\Delta t}{h_k^{n,l+1}} \left( \mathbf{F}_{k+\frac 1 2 }(\bfU_{k}^{n,l},\bfU_{k+1}^{n,l}) - \mathbf{F}_{k-\frac 1 2}(\bfU_{k-1}^{n,l},\bfU_{k}^{n,l})  \right) \right)
\]%\end{multline*}
with $h_k^{n,l+1}$ denoting the current iterate / approximation to $h_k^{n+1}$. In each iteration of the dual time stepping, this quantity is updated. 

Note that this is another advantage of using dual time stepping to solve the implicit system compared to using for example Newton method -- that we do not need to write $h_k^{n+1}$ in terms of data $\bfU^{n+1}$ and differentiate through but that we can simply reuse the same formula for updating the cell widths as for standard explicit time stepping.

%
% Start Block Time-accurate local time stepping
%
\ifthenelse{\boolean{false}}{

\subsection{Time-accurate local time stepping}\label{subsec: time-acc lts}

\todo[inline]{Den abschnitt einfach rauslassen?}

We briefly want to sketch the idea of time-accurate local time stepping for handling the presence of very small cells as we will compare in our numerical results with that approach. Full details can be found with M\"uller and Stiriba \cite{Mueller2007}. We only need a reduced version here. Again, imagine a mesh as sketched in figure \ref{fig:1d model problem} with many bigger cells and one tiny cell. The idea is to use explicit time stepping and update all cells with indices $i \le k-2$  and $i\ge k+2$ using a time step $\Delta t_{\text{ref}} \le \textcolor{red}{XXXX}$, just like in Step 1 of \textcolor{red}{XXX}. But the solution values on cells
$I_{k-1}, I_k, I_{k+1}$ are updated differently. 

\todo[inline]{SM @Ferdinand: ich verstehe nicht wie die fluesse $\mathcal{F}^\nu_{i_0\pm 1/2}$ berechnet werden $\Rightarrow$ der teil hier ist evtl falsch / muesste von dir editiert werden}

The idea is to first compute $\bfU_k^{n+1}$ and then $\bfU_{k \pm 1}^{n+1}$. To update from $\bfU_k^{n}$ to $\bfU_k^{n+1}$, a time-accurate time stepping with local step length $\Delta \tau = \alpha \Delta t_{\text{ref}}$ is used (assuming that the maximum wave speed does not change in the time interval $(t^n,t^{n+1})$). Therefore, to go from $t^n$ to $t^{n+1}$, $\frac{1}{\alpha}$ local time steps are required. Different to the local time stepping applied in dual time above, the fluxes for the update of the small cell must be computed time-accurate. This is sketched in figure \ref{fig:local_time_stepping}. \Sandra{@Ferdinand: hier bitte kurze erklaerung einfuegen.} The small cell is updated according to formula (\ref{godunov_meth:nonunif}) with the time step $\Delta\tau$
and the fluxes $\mathcal{F}^\nu_{k\pm 1/2}$ at the current time level $\tau^\nu  = t^n + \nu\Delta\tau$ with $\nu = 0,\dots,n_0$ and $\tau^{n_0} = t^{n+1}$ using 
\[
     \mathbf{U}^{n,\nu+1}_{k} =
    \frac{h_{k}^\nu}{h_{k}^{\nu+1}}\mathbf{U}^{n,\nu}_{k}
    - \frac{\Delta\tau}{h_{k}^{\nu+1}}\left[\mathcal{F}^\nu_{k+\frac{1}{2}} - \mathcal{F}^\nu_{k-\frac{1}{2}}\right].
    \label{godunov_meth:nonunif_lts}
\]
The values of the neighbouring cells remain unchanged throughout the time evolution. To ensure conservation the weighted sum of all the fluxes used in the update of cell $I_k$ need to be used to compute $\mathbf{F}_{k\pm \frac 1 2}$. This is done as follows: 
At $t^n$ the fluxes $\mathbf{F}_{k\pm 1/2}$ are computed as always. At every further time level we update the fluxes as follows
\begin{align*}
    \mathbf{F}^{\nu}_{k\pm 1/2} = \mathbf{F}^{\nu-1}_{k\pm 1/2} + \frac{\Delta\tau}{\Delta t}\mathcal{F}^{\nu}_{k\pm 1/2},\;\nu = 1,\dots,n_0.
\end{align*}
Finally, $\mathbf{F}_{k\pm \frac 1 2}^{n_0}$ together with the standard $\mathbf{F}_{k\pm \frac 3 2}$ is used to compute the update from
$\bfU_{k \pm 1}^n$ to $\bfU_{k \pm 1}^{n+1}.$

\begin{figure}[h!]
    \hspace*{0.5cm}
    \begin{tikzpicture}
        %
        %space axis at t
        \draw[color=black] (0,0) node[left] {$t^n$};
        \draw[color=black] (0,0) -- (3,0);
        \draw[color=myblue] (3,0) -- (4,0);
        \draw[color=black] (4,0) -- (7,0);
        \draw[mark=*] plot coordinates {(1,0)};
        \draw[mark=*,color=myblue] plot coordinates {(3,0)};
        \draw[mark=*,color=myblue] plot coordinates {(4,0)};
        \draw[mark=*] plot coordinates {(6,0)};
        \draw[color=black] (2,0) node[below] {$U^n_{i_0 - 1}$};
        \draw[color=myblue] (3.5,0) node[below] {$U^n_{i_0}$};
        \draw[color=black] (5,0) node[below] {$U^n_{i_0 + 1}$};
        %connect time levels
        \draw[color=black] (1,0) -- (1,2);
        \draw[color=myblue] (3,0) -- (3,2);
        \draw[color=myblue] (4,0) -- (4,2);
        \draw[color=black] (6,0) -- (6,2);
        %space axis at t + Delta tau
        \draw[color=myblue] (0,2) node[left] {$\tau^1 = t^n + \Delta\tau$};
        \draw[color=myblue] (3,2) -- (4,2);
        \draw[mark=*,color=myblue] plot coordinates {(3,2)};
        \draw[mark=*,color=myblue] plot coordinates {(4,2)};
        \draw[color=myblue] (3.5,2) node[below] {\small$U^{n,1}_{i_0}$\normalsize};
        %connect time levels
        \draw[color=black,dotted] (1,2) -- (1,3);
        \draw[color=black,dotted] (3,2.1) -- (3,3);
        \draw[color=black,dotted] (4,2.1) -- (4,3);
        \draw[color=black,dotted] (6,2) -- (6,3);
        %space axis at \tau^{n_0 - 1}
        \draw[color=myblue] (0,3) node[left] {$\tau^{n_0 - 1}$};
        \draw[color=myblue] (3,3) -- (4,3);
        \draw[mark=*,color=myblue] plot coordinates {(3,3)};
        \draw[mark=*,color=myblue] plot coordinates {(4,3)};
        \draw[color=myblue] (3.5,3) node[below] {\small$U^{n,n_0-1}_{i_0}$\normalsize};
        %connect time levels
        \draw[color=black] (1,3) -- (1,5);
        \draw[color=myblue] (3,3) -- (3,5);
        \draw[color=myblue] (4,3) -- (4,5);
        \draw[color=black] (6,3) -- (6,5);
        %space axis at t
        \draw[color=black] (0,5) node[left] {$\tau^{n_0} = t^{n + 1}$};
        \draw[color=black] (0,5) -- (7,5);
        \draw[color=myblue] (3,5) -- (4,5);
        \draw[color=black] (4,5) -- (7,5);
        \draw[mark=*] plot coordinates {(1,5)};
        \draw[mark=*,color=myblue] plot coordinates {(3,5)};
        \draw[mark=*,color=myblue] plot coordinates {(4,5)};
        \draw[mark=*] plot coordinates {(6,5)};
        \draw[color=black] (2,5) node[above] {$U^{n+1}_{i_0 - 1}$};
        \draw[color=myblue] (3.5,5) node[above] {$U^{n+1}_{i_0}$};
        \draw[color=black] (5,5) node[above] {$U^{n+1}_{i_0 + 1}$};
        %
        %draw arrows for flux computation
        %left
        \draw[color=myred] (2,0) -- (3,1);
        \draw[color=myred] (2,0) -- (3,4);
        %center-left
        \draw[color=myred] (3.5,0) -- (3,1);
        \draw[mark=triangle*,color=myred] plot coordinates {(3,1)};
        \draw[color=myred] (2.9,1) node[left] {\small$\mathcal{F}^1_{i_0 - 1/2}$\normalsize};
        \draw[color=myred] (3.5,3) -- (3,4);
        \draw[mark=triangle*,color=myred] plot coordinates {(3,4)};
        \draw[color=myred] (3,4) node[left] {\small$\mathcal{F}^{n_0-1}_{i_0 - 1/2}$\normalsize};
        %center-right
        \draw[color=myred] (3.5,0) -- (4,1);
        \draw[mark=triangle*,color=myred] plot coordinates {(4,1)};
        \draw[color=myred] (4.05,1) node[right] {\small$\mathcal{F}^1_{i_0 + 1/2}$\normalsize};
        \draw[color=myred] (3.5,3) -- (4,4);
        \draw[mark=triangle*,color=myred] plot coordinates {(4,4)};
        \draw[color=myred] (4,4) node[right] {\small$\mathcal{F}^{n_0-1}_{i_0 + 1/2}$\normalsize};
        %right
        \draw[color=myred] (5,0) -- (4,1);
        \draw[color=myred] (5,0) -- (4,4);
\end{tikzpicture}
    \caption{Simplified local time stepping for isolated small cells according to \cite{Mueller2007} \Sandra{edit fig?}}
    \label{fig:local_time_stepping}
\end{figure}

}{}
%
% END Block Time-accurate local time stepping
%

%
% Start Block Ferdinand: grid alignment and local time stepping
%
\ifthenelse{\boolean{false}}{

\textcolor{blue}{
We now want to comment on the first stated difficulty that we pointed out in the context of two phase flows with a sharp interface.
For this reason we assume that we have constant cell values $\mathbf{U}_i^n$ belonging to two different phases and that the phase boundaries are aligned with the cell boundaries.
Let us focus on a single phase boundary moving at the velocity $w \neq 0$ located at a certain $x_{i_0 + 1/2}$.
As explained in Section \ref{sec:num_probl} the phase boundary will lie inside a cell at $t^{n+1} = t^n + \Delta t$ and hence averaging will lead to unphysical phase states.
In order to avoid this situation we always align our grid to the phase boundary. This changes the size of the cells and thus we have to deal with different cell sizes.
This is done in the following way. We move the cell boundary according to the phase boundary, i.e.\ $x_{i_0 + 1/2} = x_{i_0 + 1/2} + w\Delta t$, see Figure \ref{fig:align_grid_to_pb}.}
\begin{figure}[h!]
    \hspace*{0.5cm}
    \begin{tikzpicture}
        %
        %phase boundary
        \draw[color=myred] (3,0) -- (4,2);
        %space axis at t
        \draw[color=black] (0,0) node[left] {$t^n$};
        \draw[color=black] (0,0) -- (6,0);
        \draw[mark=*] plot coordinates {(1,0)};
        \draw[mark=*] plot coordinates {(3,0)};
        \draw[mark=*] plot coordinates {(5,0)};
        \draw[color=black] (2,0) node[below] {$I_{i_0}$};
        \draw[color=black] (3,0) node[below] {\small$x_{i_0 + 1/2}$\normalsize};
        \draw[color=black] (4,0) node[below] {$I_{i_0+1}$};

        %space axis at t + Delta t
        \draw[color=black] (0,2) node[left] {$t^{n+1}$};
        \draw[color=black] (0,2) -- (6,2);
        \draw[mark=*] plot coordinates {(1,2)};
        \draw[mark=*] plot coordinates {(3,2)};
        \draw[mark=*,color=myred] plot coordinates {(4,2)};
        \draw[mark=*] plot coordinates {(5,2)};
        \draw[-latex,shorten <= 2pt,shorten >= 2pt] (3,2) arc [radius=0.7, start angle = 135, end angle = 45];
        \draw[color=black] (3.5,2.25) node[above] {$w\Delta t$};
        %aligned grid
        \draw[-latex,shorten <= 2pt,shorten >= 2pt] (-0.75,2) arc [radius=0.7, start angle = 225, end angle = 135];
        \draw[color=black] (-0.9,2.5) node[left] {align grid};
        \draw[color=black] (0,3) node[left] {$t^{n+1}$};
        \draw[color=black] (0,3) -- (6,3);
        \draw[mark=*] plot coordinates {(1,3)};
        \draw[mark=o] plot coordinates {(3,3)};
        \draw[mark=*] plot coordinates {(4,3)};
        \draw[mark=*] plot coordinates {(5,3)};
        \draw[color=black] (2.5,3) node[below] {$I_{i_0}$};
        \draw[color=black] (4,3) node[above] {\small$x_{i_0 + 1/2} + w\Delta t$\normalsize};
        \draw[color=black] (4.5,3) node[below] {$I_{i_0+1}$};
\end{tikzpicture}
    \caption{Align grid to the new position of the phase boundary.}
    \label{fig:align_grid_to_pb old}
\end{figure}
\textcolor{blue}{
Next we compare every cell to certain thresholds given by the average cell size $\Delta x_{av} := (b-a)/N$ and the two parameters $\epsilon_1$ and $\epsilon_2$.
%\footnote{For our purposes we worked with $\epsilon_1 =1.5$ and $\epsilon_2 = 0.5$.}.
Now we adjust the grid according to the following rules.
If a cell fulfills $\Delta x_i > \epsilon_1\Delta x_{av}$ the cell is too large and we split it into two.
If a cell fulfills $\Delta x_i < \epsilon_2\Delta x_{av}$ the cell is to small and we merge it with a neighbouring cell, if both belong to the same phase.
The cell values have to be updated accordingly. In the case of nucleation and cavitation we may encounter small cells that \emph{must not} be merged with the neighbouring cells.
How we deal with this situation will be explained later on.\\
\newline
Since cell sizes may change in every time step, we need to adjust equation (\ref{godunov_meth:unif}) so that it is still in conservative form.
The result reads
\begin{align}
    \mathbf{U}^{n+1}_i &= \frac{\Delta x_i^n}{\Delta x_i^{n+1}}\mathbf{U}^n_i - \frac{\Delta t}{\Delta x_i^{n+1}}\left[\mathbf{F}_{i+\frac{1}{2}} - \mathbf{F}_{i-\frac{1}{2}}\right].
    \label{godunov_meth:nonunif}
\end{align}
with $\Delta x_i^n$ and $\Delta x_i^{n+1}$ being the cell sizes at the times $t^n$ and $t^{n+1}$.
One immediately verifies that (\ref{godunov_meth:nonunif}) reduces to (\ref{godunov_meth:unif}) for constant cell sizes. Further we calculate and update the time step as follows.
The initial time step is obtained by using
\begin{align}
    S^0_{max} = \max\{|u - a|,|u + a|\} \quad\text{and}\quad \Delta t = C_{CFL}\frac{\Delta x_0}{S^0_{max}}.\label{def:init_timestep}
\end{align}
Where the maximum is computed considering all cell values. As the computation evolves we recalculate the time step as
\begin{align}
    \Delta x_{min} &= \min\limits_{i=1,\dots,N}\left\{\Delta x_i\, |\, \Delta x_i \geq \epsilon_2\Delta x_{av}\right\},\notag\\
    S^n_{max} &= \max\limits_{S}|S|,\notag\\
    \Delta t &= C_{CFL}\frac{\Delta x_{min}}{S^n_{max}}.\label{def:timestep}
\end{align}
\newline
It remains to discuss the case if a small cell of one phase is surrounded by cells of the other phase and thus must not be merged with its neighbours.
This is the case when phase creation occurs
\footnote{It is further possible that this happens when two phase boundaries move towards each other and the enclosed phase might vanish.
This situation is not considered here and will be discussed in the future.}.
In this case we proceed as follows.
For cells with size $\Delta x_i < \epsilon_2\Delta x_{av}$ we apply a \emph{local time stepping method} as properly discussed by M\"uller and Stiriba \cite{Mueller2007}.}
\begin{figure}[h!]
    \hspace*{0.5cm}
    \begin{tikzpicture}
        %
        %space axis at t
        \draw[color=black] (0,0) node[left] {$t^n$};
        \draw[color=black] (0,0) -- (3,0);
        \draw[color=myblue] (3,0) -- (4,0);
        \draw[color=black] (4,0) -- (7,0);
        \draw[mark=*] plot coordinates {(1,0)};
        \draw[mark=*,color=myblue] plot coordinates {(3,0)};
        \draw[mark=*,color=myblue] plot coordinates {(4,0)};
        \draw[mark=*] plot coordinates {(6,0)};
        \draw[color=black] (2,0) node[below] {$U^n_{i_0 - 1}$};
        \draw[color=myblue] (3.5,0) node[below] {$U^n_{i_0}$};
        \draw[color=black] (5,0) node[below] {$U^n_{i_0 + 1}$};
        %connect time levels
        \draw[color=black] (1,0) -- (1,2);
        \draw[color=myblue] (3,0) -- (3,2);
        \draw[color=myblue] (4,0) -- (4,2);
        \draw[color=black] (6,0) -- (6,2);
        %space axis at t + Delta tau
        \draw[color=myblue] (0,2) node[left] {$\tau^1 = t^n + \Delta\tau$};
        \draw[color=myblue] (3,2) -- (4,2);
        \draw[mark=*,color=myblue] plot coordinates {(3,2)};
        \draw[mark=*,color=myblue] plot coordinates {(4,2)};
        \draw[color=myblue] (3.5,2) node[below] {\small$U^{n,1}_{i_0}$\normalsize};
        %connect time levels
        \draw[color=black,dotted] (1,2) -- (1,3);
        \draw[color=black,dotted] (3,2.1) -- (3,3);
        \draw[color=black,dotted] (4,2.1) -- (4,3);
        \draw[color=black,dotted] (6,2) -- (6,3);
        %space axis at \tau^{n_0 - 1}
        \draw[color=myblue] (0,3) node[left] {$\tau^{n_0 - 1}$};
        \draw[color=myblue] (3,3) -- (4,3);
        \draw[mark=*,color=myblue] plot coordinates {(3,3)};
        \draw[mark=*,color=myblue] plot coordinates {(4,3)};
        \draw[color=myblue] (3.5,3) node[below] {\small$U^{n,n_0-1}_{i_0}$\normalsize};
        %connect time levels
        \draw[color=black] (1,3) -- (1,5);
        \draw[color=myblue] (3,3) -- (3,5);
        \draw[color=myblue] (4,3) -- (4,5);
        \draw[color=black] (6,3) -- (6,5);
        %space axis at t
        \draw[color=black] (0,5) node[left] {$\tau^{n_0} = t^{n + 1}$};
        \draw[color=black] (0,5) -- (7,5);
        \draw[color=myblue] (3,5) -- (4,5);
        \draw[color=black] (4,5) -- (7,5);
        \draw[mark=*] plot coordinates {(1,5)};
        \draw[mark=*,color=myblue] plot coordinates {(3,5)};
        \draw[mark=*,color=myblue] plot coordinates {(4,5)};
        \draw[mark=*] plot coordinates {(6,5)};
        \draw[color=black] (2,5) node[above] {$U^{n+1}_{i_0 - 1}$};
        \draw[color=myblue] (3.5,5) node[above] {$U^{n+1}_{i_0}$};
        \draw[color=black] (5,5) node[above] {$U^{n+1}_{i_0 + 1}$};
        %
        %draw arrows for flux computation
        %left
        \draw[color=myred] (2,0) -- (3,1);
        \draw[color=myred] (2,0) -- (3,4);
        %center-left
        \draw[color=myred] (3.5,0) -- (3,1);
        \draw[mark=triangle*,color=myred] plot coordinates {(3,1)};
        \draw[color=myred] (2.9,1) node[left] {\small$\mathcal{F}^1_{i_0 - 1/2}$\normalsize};
        \draw[color=myred] (3.5,3) -- (3,4);
        \draw[mark=triangle*,color=myred] plot coordinates {(3,4)};
        \draw[color=myred] (3,4) node[left] {\small$\mathcal{F}^{n_0-1}_{i_0 - 1/2}$\normalsize};
        %center-right
        \draw[color=myred] (3.5,0) -- (4,1);
        \draw[mark=triangle*,color=myred] plot coordinates {(4,1)};
        \draw[color=myred] (4.05,1) node[right] {\small$\mathcal{F}^1_{i_0 + 1/2}$\normalsize};
        \draw[color=myred] (3.5,3) -- (4,4);
        \draw[mark=triangle*,color=myred] plot coordinates {(4,4)};
        \draw[color=myred] (4,4) node[right] {\small$\mathcal{F}^{n_0-1}_{i_0 + 1/2}$\normalsize};
        %right
        \draw[color=myred] (5,0) -- (4,1);
        \draw[color=myred] (5,0) -- (4,4);
\end{tikzpicture}
    \caption{Simplified local time stepping for isolated small cells according to \cite{Mueller2007}}
    \label{fig:local_time_stepping}
\end{figure}
\textcolor{blue}{
We only need a simple reduced version of the method presented in \cite{Mueller2007} and will briefly present the used method in the following, see also Figure \ref{fig:local_time_stepping}.\\
Let $i_0$ be the index of the small cell. We then consider the triplet $\{U_{i_0-1},U_{i_0},U_{i_0+1}\}$ and the corresponding quantities.
Now we perform a time evolution starting at $t^n$ with the small time steps $\Delta\tau$ until the final time $t^{n+1} = t^n + \Delta t$.
The small cell is updated according to formula (\ref{godunov_meth:nonunif}) with the time step $\Delta\tau$
and the fluxes $\mathcal{F}^\nu_{i_0\pm 1/2}$ at the current time level $\tau^\nu  = t^n + \nu\Delta\tau$ with $\nu = 0,\dots,n_0$ and $\tau^{n_0} = t^{n+1}$
\footnote{ A common choice is to choose $\Delta \tau = 2^{-n_0}\Delta t$ such that the \emph{CFL-condition} on the fine level is satisfied, see \cite{Mueller2007}.}. We have
\begin{align}
    \mathbf{U}^{n,\nu+1}_{i_0} &=
    \frac{\Delta x_{i_0}^\nu}{\Delta x_{i_0}^{\nu+1}}\mathbf{U}^{n,\nu}_{i_0}
    - \frac{\Delta\tau}{\Delta x_{i_0}^{\nu+1}}\left[\mathcal{F}^\nu_{i_0+\frac{1}{2}} - \mathcal{F}^\nu_{i_0-\frac{1}{2}}\right].
    \label{godunov_meth:nonunif_lts}
\end{align}
The values of the neighbouring cells remain unchanged throughout the time evolution.
For the update of these (large) neighbouring cells we have to take care of the fluxes at $x_{i_0-1/2}$ and $x_{i_0+1/2}$.
According to \cite{Mueller2007} this is done in the following way. At $t^n$ we calculate the fluxes $\mathbf{F}_{i_0\pm 1/2}$ as usual. At every further time level we update these fluxes as follows
\begin{align*}
    \mathbf{F}^{\nu}_{i_0\pm 1/2} = \mathbf{F}^{\nu-1}_{i_0\pm 1/2} + \frac{\Delta\tau}{\Delta t}\mathcal{F}^{\nu}_{i_0\pm 1/2},\;\nu = 1,\dots,n_0.
\end{align*}
The fluxes obtained at the end of the evolution of the small cell may then be used in (\ref{godunov_meth:nonunif}) to advance the neighbouring cells from $t^n$ to $t^{n+1}$.}

%
% END Block Ferdinand: grid alignment and local time stepping
%
}{}

%% file: examples.tex
\section{Numerical results}\label{sec:num_ex}

In this section we present various results using the mixed explicit implicit approach described in section \ref{sec:num_interface} together with dual time stepping to solve the resulting implicit systems. We will first consider a single phase flow problem, where we discuss our choice of parameters in the dual time stepping. %Then we will present a two phase problem without phase transition.
Afterwards, we will present tests involving nucleation and cavitation. Here, the focus is on verifying that the results produced with our new approach are physically meaningful and that the new approach is relatively fast. 

All results for two phase flow problems are compared to the solution of the Riemann problem obtained with the Newton method using the tolerance $\varepsilon_{tol} = 10^{-9}$.
More precisely, since in some cases the standard Newton method failed to converge, we used the \emph{Newton-Armijo} method as presented in Kelley \cite{Kelley1995}.
The saturation pressure is calculated using the steam tables Ref.\ \cite{Wagner1998}.
If not stated otherwise we used the exact solution of the Riemann problem at the phase boundary in our numerical calculations as Riemann solvers. In the two phase flow examples away from the phase boundary any Riemann solver can be used. Here we applied for the numerical solution of the two phase flow examples the HLL solver as presented in Toro\cite{Toro2009}. In our first test, where the full compressible Euler equations are solved, we use an exact Riemann solver.

For the isothermal case we model the vapor phase as an ideal gas using the relation
\begin{align}
    p_V = \frac{kT_0}{m}\rho_V,\quad m = \frac{2\cdot1.0079 + 15.9994}{6.02205\cdot10^{26}}\,\si{kg},\quad k = 1.380658\cdot10^{-23}\,\si{JK^{-1}},\label{eos:isoT_idgas}
\end{align}
where $k$ is the Boltzmann constant and $m$ the mass of a single water molecule.
The parameter $\tau$ in the kinetic relation \eqref{isoT:kin_rel1} is given by
\begin{align*}
    \tau = \frac{1}{\sqrt{2\pi}}\left(\frac{m}{kT_0}\right)^{3/2}.
\end{align*}
The liquid phase is described using the linear Tait EOS
\begin{align}
    p_L = p_0 + K_0\left(\frac{\rho_L}{\rho_0} - p_0\right).\label{eos:isoT_lintait}
\end{align}
All quantities with a subscript zero are saturation quantities at the given temperature $T_0$ evaluated according to the steam tables in Ref. \cite{Wagner1998}.
The EOS for the individual phases are then connected linearly with respect to the density and it is assured that the Maxwell condition holds.

Concerning the parameters in our dual time stepping approach, we use $\theta_{k \pm 1} = 0.9$ and $\theta_k = \alpha$ for setting the local time step lengths on the small cell and its two neighbors, compare also algorithm 2 in table \ref{tab: algorithms}. For the \textit{stopping criterium} in the dual time stepping, we define the residuum on cell $I_i$ as (with the division to be understand in a component-wise way)
\[
\text{res}_i = 
\left \lVert \frac{\left(1 +  \frac{\Delta \tau_i}{\Delta t}\right) \left( \bfU_{i}^{n,l+1} - \bfU_{i}^{n,l} \right)}{\Delta \tau_i \max(\bfU_i^{n,l},1)}\right\rVert_{\infty}.
\]
Note the division by the time step length $\Delta \tau_i$ of the local dual time stepping. As our updates are designed to take a time step of that length, it is important to scale that out by checking for steady state. We also divide by $\max(\bfU_i^{n,l},1)$ to account for large differences (in terms of absolute sizes) of density and momentum when comparing the residuals of cells that belong to different phases. Finally, we require as stopping criteria $\max(\text{res}_{k-1},\text{res}_{k+1}) < 10^{-3}$ and $\text{res}_{k} < 10^{-1}$.  

We found that this works well in our numerical tests. These stopping tolerances seem very imprecise. But due to dividing by $\Delta \tau_{i}$, the residual goes up to $10^8$ or higher at initial iterations if $\alpha \ll 1$. We also note that one needs to be careful to stay away from computing on numerical artefacts/round-off errors. For the nucleation test in section \ref{subsubsec:ex3} below, we get 
$\Delta \tau_k = 6 \cdot 10^{-12}$ at initial iterations. Assuming that $\bfU_k^{n,l} = \mathcal{O}(1)$, there are only few meaningful digits left in the double precision arithmetic.

\subsection{Stationary Phase Boundary without Phase Transition}\label{subsec: Euler test}
In our first test, we consider a single phase problem. The goal of this test is to examine 
 how the choice of the parameters $\theta_i, i=k{-1},k,k{+}1,$ influences the number of iterations in the dual time stepping. 
 By means of $\Delta \tau_i = \theta_i \Delta t_{\text{ref}}$, these parameters control the local time step lengths in the dual time stepping, compare also algorithm 2 in table \ref{tab: algorithms}.

We consider the full one-dimensional compressible Euler equations given by
\[
\begin{pmatrix}
    \rho \\ \rho u \\ E
    \end{pmatrix}_t + 
    \begin{pmatrix} \rho u \\ \rho u^2 + p \\ (E + p)u \end{pmatrix}_x = 
    \begin{pmatrix}
    0 \\ 0 \\0
    \end{pmatrix}
\] 
on a fixed mesh. We solve the Lax shock tube test for ideal gas using $\gamma=1.4$ given by the initial conditions
\[
(\rho,u,p)(x,0) = 
\begin{cases}
(0.445, 0.698, 3.528), & x < 0,\\
(0.5, 0, 0.571), & x > 0,
\end{cases}
\]
% DL = .445, DR = .5
    % UL = .698, UR = 0
    % PL = 3.528, PR = .571
on the domain $[x_L,x_R] = [-5,5]$ for the model mesh shown in figure \ref{fig:1d model problem}. We place the left boundary of the small cell $I_k$ at $x=0$, i.e., $x_{k-\frac 1 2}=0$. We define $h = (x_R-x_L)/N$ with $N$ being the number of equidistant cells that we use in our computations. (This setup results in actually using the domain $[-5,5+\alpha h]$ in our tests.) We run 
until time $t_{\text{end}}=1.0$ using $C_{CFL}=0.8$.
%
% \begin{table}[h!]
%     %
%     \centering
%     \renewcommand*{\arraystretch}{1.25}
%     %
%     \begin{tabular}{cc|ccc}
%         %
%          $\theta_{k \pm 1}$ & $N$         & $\alpha=10^{-2}$          & $\alpha=10^{-4}$              & $\alpha=10^{-6}$          \\
%         \hline
%         0.9   &   $N=100$    & 23 & 42 & 63 \\
%               &   $N=200$    & 18 & 36 & 54 \\
%               &   $N=400$    & 16 & 32 & 49 \\
%               \hline
%         0.1   &   $N=100$    & 74 & 82 & 119 \\
%               &   $N=200$    & 65 & 73 & 109 \\
%               &   $N=400$    & 54 & 64 & 100 \\    
%               \hline
%         $\alpha$  & $N=100$  & 681 & 67\,607 & * \\
%               &   $N=200$    & 595 & 59\,029 & * \\
%               &   $N=400$    & 492 & 48\,857 & * \\                      
%         %
%     \end{tabular}
%     %
%     \caption{Lax shock tube: Average number of iterations in dual time stepping for $\theta_k = \alpha$ and varying choices for % $\theta_{k \pm 1}$. (The results marked with * took too long to compute in a reasonable time.) --- \Sandra{OLD}}
%     \label{tab:iterations}
%     %
% \end{table}
%
\begin{table}[h!]
    \centering
    \renewcommand*{\arraystretch}{1.25}
    \begin{tabular}{cc|ccc}
         $\theta_{k \pm 1}$ & $N$         & $\alpha=10^{-2}$          & $\alpha=10^{-4}$              & $\alpha=10^{-6}$          \\
        \hline
        0.9   &   $N=100$    & 16 & 35 & 57 \\
              &   $N=200$    & 14 & 32 & 51 \\
              &   $N=400$    & 13 & 30 & 47 \\
              \hline
        0.1   &   $N=100$    & 48 & 59 & 97 \\
              &   $N=200$    & 38 & 51 & 90 \\
              &   $N=400$    & 29 & 44 & 81 \\    
              \hline
        $\alpha$  & $N=100$  & 439 & 43\,470 & * \\
              &   $N=200$    & 347 & 34\,356 & * \\
              &   $N=400$    & 246 & 24\,178 & * \\                      
    \end{tabular}
    \caption{Lax shock tube: Average number of iterations in dual time stepping for $\theta_k = \alpha$ and varying choices for $\theta_{k \pm 1}$. (The results marked with * took too long to compute in a reasonable time.)}
    \label{tab:iterations}
\end{table}
In table \ref{tab:iterations}, we show the number of iterations needed in the dual time stepping, averaged over the running time $[0,T]$, for different setups. Due to stability, we need to choose $\theta_k \in (0,\alpha].$ We use the biggest possible time step $\theta_k = \alpha$ here. For $\theta_{k \pm 1}$ there holds for stability reasons $\theta_{k \pm 1} \in (0,1]$, i.e., there is some freedom of choice here. Comparing the results for $\theta_{k \pm 1} = 0.9$ and $\theta_{k \pm 1} = 0.1$, we observe that it is better to use bigger values for $\theta_{k \pm 1}$. We made a similar observation in other tests. 

To contrast, we also present results for choosing $\theta_{k \pm 1} = \alpha$, corresponding to advancing all three cells $I_{k-1},I_k,I_{k+1}$ simultaneously based on the allowed CFL number of cell $I_k$. Here, it is clearly visible that the number of iterations scales roughly with $\mathcal{O}(\frac{1}{\alpha})$. For $\alpha = 10^{-5}$, we get %676\,027, 590\,248, and 488\,532 
434\,663, 343\,526, and 241\,731
average iterations for $N=100$, $N=200$, and $N=400$, respectively. Note that we do not observe the scaling with $\mathcal{O}(\frac{1}{\alpha})$ for using proper local time stepping with different time step lengths, i.e., when setting $\theta_{k \pm 1} = 0.9$ for example. We only observe a mild increase in the number of iterations. This is essential as the time-accurate local time stepping (when using explicit time stepping everywhere) naturally requires roughly $\mathcal{O}(\frac{1}{\alpha})$ iterations for doing one \textit{big} time step of length $\Delta t_{\text{ref}}$ on the small cell $I_k$, and this is exactly the behavior that we want to avoid.

\subsection{Cavitation}\label{subsubsec:ex4}
Next we discuss the case of cavitation. Consider the initial data given in Table \ref{tab:init_dat_4}. % and it corresponds to \emph{Example 4} in \cite{Hantke2013}.
\begin{table}[h!]
    \centering
    \renewcommand*{\arraystretch}{1.25}
    \begin{tabular}{c|cccc}
                     & $p_L^-$          & $u_L^-$                 & $p_L^+$          & $u_L^+$ \\%                & $T_0$            & $p_0$\\
        \hline
        Initial Data & $60000\,\si{Pa}$ & $-4\,\si{\frac{m}{s}}$ & $60000\,\si{Pa}$ & $4\,\si{\frac{m}{s}}$ %& $363.15\,\si{K}$ & $70182.360745\,\si{Pa}$
    \end{tabular}
    \begin{tabular}{c|cc}
                          & $T_0$            & $p_0$\\
        \hline
        Saturation Values & $363.15\,\si{K}$ & $70182.360745\,\si{Pa}$
    \end{tabular}
    \caption{Cavitation Test: Initial Data.}
    \label{tab:init_dat_4}
\end{table}

The following computation was performed with $C_{CFL} = 0.5,\;h^0 = 10^{-2}\,\si{m},\; x \in [-2,2]\,\si{m}$, and $t_{\text{end}} = 5\cdot10^{-4}\,\si{s}$.
The time step is calculated according to the CFL-condition \eqref{cfl_cond:unif} and the phases are distinguished using the values given in Table \ref{tab:max_vap_lintait363}. 
\begin{table}[h!]
    \centering
    \renewcommand*{\arraystretch}{1.5}
    \begin{tabular}{cc||cc}
        $\tilde{p}$             & $\tilde{\rho}$                  & $p_{min}$    & $\rho_{min}$ \\
        \hline
        $70388.660656\,\si{Pa}$ & $0.419977\,\si{\frac{kg}{m^3}}$ & $0\,\si{Pa}$ & $965.289008\,\si{\frac{kg}{m^3}}$
    \end{tabular}
    \caption{Maximum vapor pressure and corresponding quantities for the linear Tait EOS together with the ideal gas EOS at $T_0 = 363.15\,\si{K}$.}
    \label{tab:max_vap_lintait363}
\end{table}

\begin{figure}[h!]
    \includegraphics[scale=1]{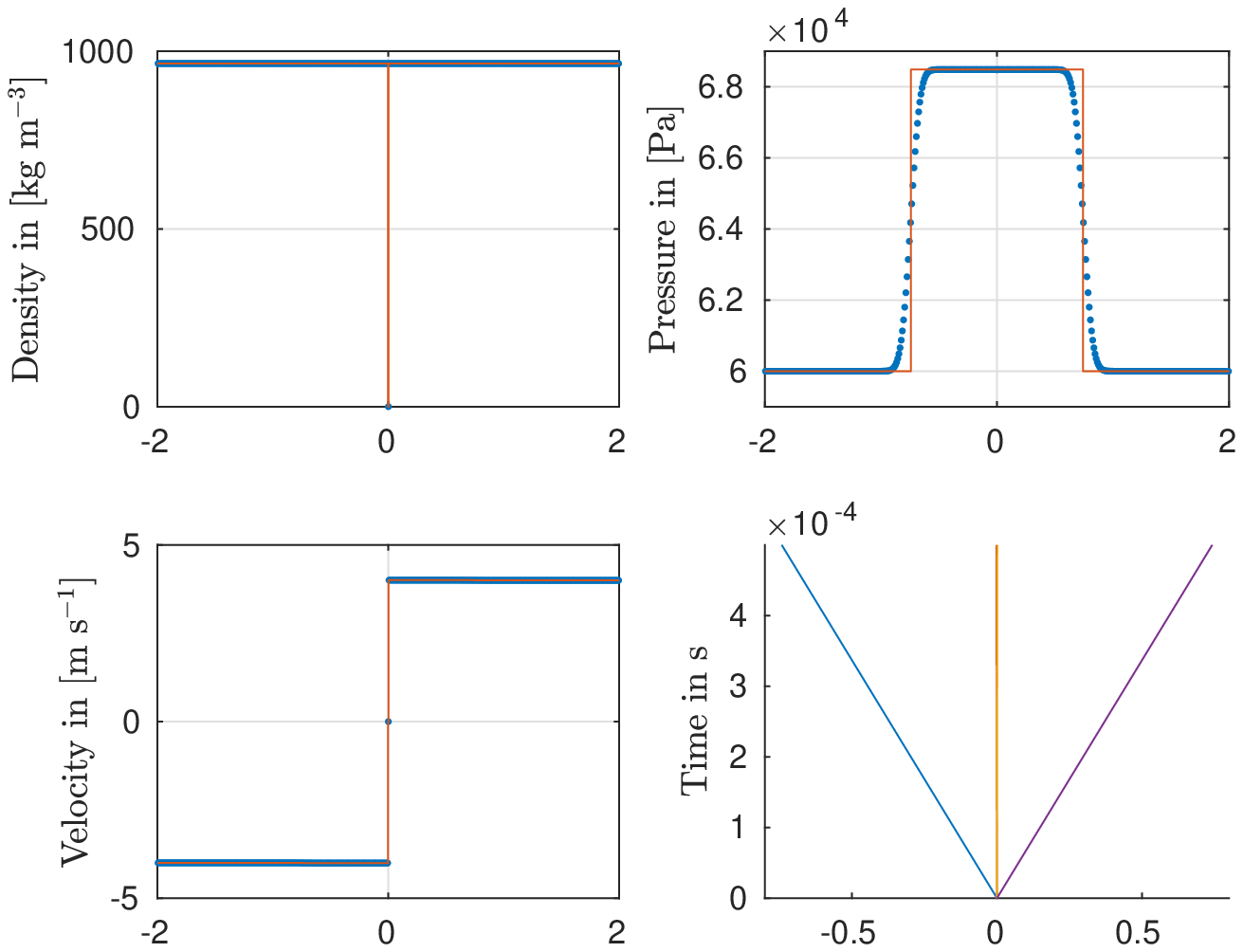}
    \caption{Cavitation Test: Computed (blue) and exact (red) solution as well as the exact wave structure of the solution showing the outgoing shock waves and the two phase boundaries in the middle of the fan. (Note that the two phase boundaries (shown in orange) are so close that they cannot be distinguished in the given picture and appear as one line.)}
  \label{fig:cav_ex}
\end{figure}

To treat the cavitation case we apply the mixed explicit implicit time stepping described above (in combination with dual time stepping).
We will also compare with using the fully explicit time-accurate local time stepping (LTS) as presented in Ref. \cite{Thein2018} based on Ref. \cite{Mueller2007}.
The results at the final time for the computed solution using the mixed explicit implicit approach together with the exact solution are given in Figure \ref{fig:cav_ex}.
For this test, the cell size of the vapor phase after its creation is $h_V = 2.701184\cdot 10^{-5}\,\si{m}$ which leads to an initial $\alpha \sim \mathcal{O}(10^{-3})$. 
The phase boundaries move at $w_{\text{left/right}} = \pm 4.007636\,\si{\frac{m}{s}}$ and thus the size of the liquid cell at time $t_{\text{end}}$ is $(w_{right} - w_{left})\cdot t_{\text{end}} = 4.007635\cdot 10^{-3}\,\si{m}.$

%\Sandra{english im folgenden satz sowie in dem davor}\Ferdinand{Lass uns drüber reden...} Note for the wave structure, that the phase boundaries are so close that they cannot be distinguished in the given picture.  \Sandra{ich verstehe diesen absatz nicht wirklich}\Ferdinand{Lass uns drüber reden...}

The values in the star region for the pressure are given in Table \ref{tab:sol_dat_4_p}. 
\begin{table}[h!]
    \centering
    \renewcommand*{\arraystretch}{1.25}
    \begin{tabular}{c|ccc}
                  & $p_L^{\ast,-}$          & $p_V^\ast$              & $p_L^{\ast,+}$          \\
        \hline
        exact     & $68483.741137\,\si{Pa}$ & $68477.181783\,\si{Pa}$ & $68483.741137\,\si{Pa}$ \\
        explicit  & $68483.741813\,\si{Pa}$ & $68477.181783\,\si{Pa}$ & $68483.740209\,\si{Pa}$ \\
        implicit  & $68483.740426\,\si{Pa}$ & $68477.181783\,\si{Pa}$ & $68483.740588\,\si{Pa}$
    \end{tabular}
    \caption{Cavitation Test: Exact and computed solution (pressure).}
    \label{tab:sol_dat_4_p}
\end{table}
\begin{table}[h!]
    \centering
    \renewcommand*{\arraystretch}{1.25}
    \begin{tabular}{c|ccc}
                  & $u_L^{\ast,-}$                 & $u_V^\ast$                                & $u_L^{\ast,+}$                \\
        \hline
        exact     & $-4.005940\,\si{\frac{m}{s}}$ & $0\,\si{\frac{m}{s}}$                     & $4.005940\,\si{\frac{m}{s}}$ \\
        explicit & $-4.005940\,\si{\frac{m}{s}}$ & $-2.127421\cdot 10^{-10}\,\si{\frac{m}{s}}$ & $4.005940\,\si{\frac{m}{s}}$ \\
        implicit & $-4.005940\,\si{\frac{m}{s}}$ & $\phantom{-}1.575386\cdot 10^{-10}\,\si{\frac{m}{s}}$ & $4.005930\,\si{\frac{m}{s}}$
    \end{tabular}
    \caption{Cavitation Test: Exact and computed solution (velocity).}
    \label{tab:sol_dat_4_u}
    %
   % \vspace*{-20pt}
    %
\end{table}
We present exact solution values obtained solving the Riemann problem as presented in appendix \ref{sec:app}, the computed results using our new approach (\textit{implicit}), as well as using explicit LTS. In Table \ref{tab:sol_dat_4_u}, we present the corresponding solution data for the 
velocity in the star region. For both tables, the values produced by these two approaches show no significant differences.

We now compare the performance of the mixed explicit implicit scheme based on dual time stepping (DTS) with using a fully explicit LTS. 
In Table \ref{tab:comp_costs_2}, we present the total number of local iterations performed during the complete computation with $166$ large time steps. %The precise behaviour of the local iterations during the course of the simulation is visualized in figure \ref{fig:loc_it_cav_ex}.
Overall, the explicit LTS needs roughly speaking a factor of 5 more iterations than the implicit approach does. This is a good speed-up considering that the size of the small cell at creation was only about a factor of $\mathcal{O}(10^{-3})$ smaller than the other cells and grew pretty quickly. The different factor in terms of iterations needed is only partially reflected in the running times as the local time stepping only accounts for $\sim 30 \%$ of the running time in the implicit setting and for $\sim 60 \%$ in the explicit setting. 

\begin{table}[h!]
    \centering
    \renewcommand*{\arraystretch}{1.25}
    \begin{tabular}{c|cccc}
                     & explicit LTS          & implicit (DTS)\\
        \hline
        number of total local iterations & $1941$ & $355$ \\
        computation time & $4.93\si{s}$ & $2.94\si{s}$
    \end{tabular}
    \caption{Cavitation Test: Comparison of computational costs.}
    \label{tab:comp_costs_2}
\end{table}
%
%
% \begin{figure}[h!]
%     %
%     \includegraphics[scale=0.6]{graphics/prob14b_loc_it_comp.eps}
%     \caption{Cavitation Test: Local iterations per time step: explicit (blue) vs.\ implicit (red).}
%   \label{fig:loc_it_cav_ex}
%     %
% \end{figure}

% \Ferdinand{Der folgende Absatz kann evtl. auch raus...}
% Note that cavitation may be a very fast process in terms of the interface velocity depending on the initial velocity of the liquid. Hence the % vapor cell might be growing rather fast and in that case the explicit local time stepping might be equally suited and computationally efficient, % cf.\ Thein\cite{Thein2018}.

\subsection{Nucleation}\label{subsubsec:ex3}

The final example is a nucleation test case with the initial data as given in Table \ref{tab:init_dat_3}. This corresponds to \emph{Example 3} in Hantke et al.\cite{Hantke2013}.
\begin{table}[h!]
    \centering
    \renewcommand*{\arraystretch}{1.25}
    \begin{tabular}{c|cccc}
                     & $p_V^-$          & $u_V^-$                 & $p_V^+$          & $u_V^+$ \\%                  & $T_0$            & $p_0$\\
        \hline
        Initial Data & $70000\,\si{Pa}$ & $2.7\,\si{\frac{m}{s}}$ & $70000\,\si{Pa}$ & $-2.7\,\si{\frac{m}{s}}$ %& $363.15\,\si{K}$ & $70182.360745\,\si{Pa}$
    \end{tabular}
    \begin{tabular}{c|cc}
                          & $T_0$            & $p_0$\\
        \hline
        Saturation Values & $363.15\,\si{K}$ & $70182.360745\,\si{Pa}$
    \end{tabular}
    \caption{Nucleation Test: Initial Data.}
    \label{tab:init_dat_3}
\end{table}

The computation was performed with $C_{CFL} = 0.5,\;h^0 = 10^{-2}\,\si{m},\; x \in [-2,2]\,\si{m}$, and $t_{\text{end}} = 5\cdot10^{-4}\,\si{s}$.
The time step is calculated according to the CFL-condition \eqref{cfl_cond:unif}. %s given in (\ref{def:init_timestep}) and (\ref{def:timestep}).
The phases are distinguished using the values given in Table \ref{tab:max_vap_lintait363}.

\begin{figure}[h!]
    \includegraphics[scale=1]{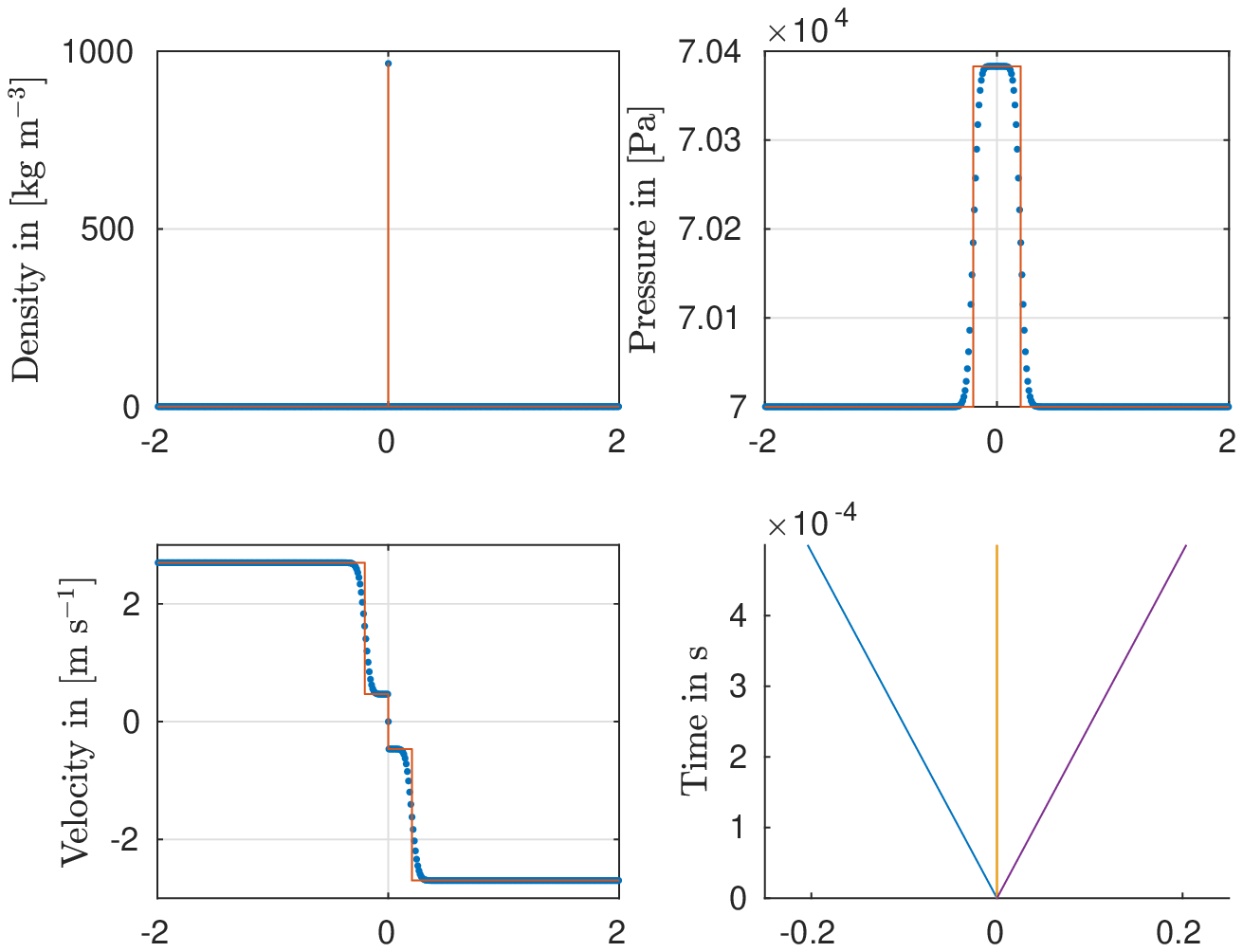}
    \caption{Nucleation Test: Computed (blue) and exact (red) solution as well as the exact wave structure of the solution showing the outgoing shock waves and the two phase boundaries in the middle of the fan. (Note that the two phase boundaries (shown in orange) are so close that they cannot be distinguished in the given picture and appear as one line.)}
  \label{fig:nucl_ex}
\end{figure}

The results for the mixed explicit implicit time stepping together with the exact solution are shown in Figure \ref{fig:nucl_ex}. 

The values in the star region for the pressure are given in Table \ref{tab:sol_dat_3_p}. We again present the exact solution as well as the solutions computed with the new explicit implicit scheme and with using explicit LTS. 
\begin{table}[h!]
    \centering
    \renewcommand*{\arraystretch}{1.25}
    \begin{tabular}{c|ccc}
                  & $p_V^{\ast,-}$          & $p_L^\ast$              & $p_V^{\ast,+}$          \\
        \hline
        exact     & $70383.024449\,\si{Pa}$ & $70383.115685\,\si{Pa}$ & $70383.024449\,\si{Pa}$ \\
        explicit  & $70382.992562\,\si{Pa}$ & $70383.083927\,\si{Pa}$ & $70382.992562\,\si{Pa}$ \\
        implicit  & $70382.992468\,\si{Pa}$ & $70383.072633\,\si{Pa}$ & $70382.992468\,\si{Pa}$
    \end{tabular}
    \caption{Nucleation Test: Exact and computed solution (pressure).}
    \label{tab:sol_dat_3_p}
\end{table}
The corresponding values for the velocity in the star region are given in Table \ref{tab:sol_dat_3_u}.
\begin{table}[h!]
    \centering
    \renewcommand*{\arraystretch}{1.25}
    \begin{tabular}{c|ccc}
                  & $u_V^{\ast,-}$               & $u_L^\ast$                                 & $u_V^{\ast,+}$                \\
        \hline
        exact     & $0.466005\,\si{\frac{m}{s}}$ & $0\,\si{\frac{m}{s}}$                      & $-0.466005\,\si{\frac{m}{s}}$ \\
        explicit  & $0.465933\,\si{\frac{m}{s}}$ & $3.161887\cdot 10^{-12}\,\si{\frac{m}{s}}$ & $-0.465933\,\si{\frac{m}{s}}$\\
        implicit  & $0.465933\,\si{\frac{m}{s}}$ & $1.2214509\cdot 10^{-10}\,\si{\frac{m}{s}}$ & $-0.465933\,\si{\frac{m}{s}}$
    \end{tabular}
    \caption{Nucleation Test: Exact and computed solution (velocity).}
    \label{tab:sol_dat_3_u}
\end{table}
Both methods show no significant differences for the computed values and are also in good agreement with the exact solution.

For this test, the cell size of the liquid phase after its creation is $h_L = 4.429422\cdot 10^{-9}\,\si{m}$ which leads to an initial $\alpha \sim \mathcal{O}(10^{-7})$. The values for the interface velocities are $w_{\text{left/right}} = \pm 0.000203\,\si{\frac{m}{s}}$ and thus the size of the liquid cell at time $t_{\text{end}}$ is $(w_{right} - w_{left})\cdot t_{\text{end}} = 2.03\cdot 10^{-7}\,\si{m}.$

In this setup, we expect bigger time differences between the implicit DTS and the explicit LTS approach than for cavitation as the run time should be dominated by the local iterations in the neighborhood of the small cell. Indeed, when comparing the total number of local iterations performed during the complete computation with $143$ large time steps, we have roughly 
a factor of 20 more iterations for the explicit approach, compare the data in table \ref{tab:comp_costs_3}.
\begin{table}[h!]
    \centering
    \renewcommand*{\arraystretch}{1.25}
    \begin{tabular}{c|cccc}
                     & explicits LTS          & implicit (DTS)\\
        \hline
        number of total local iterations & $27'154'016$ & $1'367'460$ \\
        computation time & $9811\si{s}$ & $509\si{s}$
    \end{tabular}
    \caption{Nucleation Test: Comparison of computational costs.}
    \label{tab:comp_costs_3}
\end{table}
This difference in iterations is also reflected in the running time. We now see roughly a factor of 20 there as well.

% The precise behaviour of the local iterations is visualized in figure \ref{fig:loc_it_nucl_ex}.
%
% \begin{figure}[h!]
%     %
%     \includegraphics[scale=0.6]{graphics/prob07_loc_it_comp.eps}
%     \caption{Local iterations per time step: explicit (blue) vs.\ implicit (red) - Nucleation Test \Sandra{ich wuerde die fig rauslassen -- % wirft viele Fragen auf...}\Ferdinand{Dann sollten wir das im ersten Bsp. aber auch tun.}}
%   \label{fig:loc_it_nucl_ex}
%     %
% \end{figure}
%

%% file: conclusion.tex
\section{Conclusion}\label{sec:conclusion}
In this work we suggest a new numerical method to treat nucleation and cavitation or more generally two-phase flow problems with sharp interfaces that result in the existence of tiny cells. We treat the neighborhood of the tiny cells implicitly for stability while using an explicit time stepping scheme everywhere else to keep the cost low. Using flux bounding to couple the schemes as introduced in Ref. \cite{May_Berger_explimpl} results in the mixed explicit implicit scheme being conservative and stable. We suggest to use dual time stepping for solving the resulting implicit systems in the neighborhood of the tiny cells. 
Our numerical results, which include nucleation and cavitation tests, show an accurate and robust performance of the mixed scheme. In terms of running times we observed significant speed-ups compared to a fully explicit local time stepping approach.
In the future, we plan to explore the speed-up options further by looking into better starting points for the dual time stepping (instead of simply using $\bfU^n$) or Newton-like approaches. We also plan to extend the approach to two dimensions. 

%In this work we presented and briefly discussed a numerical method to solve Riemann problems for compressible isothermal Euler equations for two phase flows with and without phase transition.
%The analysis for the exact solution was done by Hantke et al. in \cite{Hantke2013}. We highlighted major difficulties and proposed solutions to overcome these difficulties.
%Future work will include convergence and stability studies. Furthermore an extension to mixtures, nucleation and cavitation will be of interest.
%Another crucial point will be the extension to higher dimensions and non isothermal problems.

%% file: appendix.tex
\section{Solution at the Interface}\label{sec:app}
In the following we want to discuss the solution at the phase boundary between the vapor and the liquid phases. Here we rely on the results presented in Refs. \cite{Thein2018,Hantke2019a}.
The solution of the two phase Riemann problem consists of three waves, separating four constant states.
In particular we search the values in the star region, i.e.\ $\rho_V^\ast,p_V^\ast,u_V^\ast$ in the vapor phase,
$\rho_L^\ast,p_L^\ast,u_L^\ast$ in the liquid phase, and the velocity $w$ of the phase boundary or the mass flux $z$.
The pressure and the density inside each phase are not independent from each other due to the EOS given for the particular phase.
Further we know that the pressures at the interface are uniquely linked using equation
$ %$\begin{equation}
    \dbl p\dbr = -z^2\dbl v \dbr.\label{isoT:pb_mom_balance}
$ %\end{equation}
The mass flux is given by the kinetic relation (\ref{isoT:kin_rel1})
\begin{align*}
    z = \tau p_V\dbl g + e^{kin}\dbr = \tau p_V\left[g_L - g_V + e^{kin}_L - e^{kin}_V\right].
\end{align*}
Thus we have to solve the following system for $(p_V^\ast,p_L^\ast)$ to obtain the complete solution
\begin{align}
    \mathbf{0} = \mathbf{G}(p_V^\ast,p_L^\ast) = \left(\begin{matrix}
        & \dbl p\dbr + z^2\dbl v \dbr\\
        & f_V(p_V^\ast,\mathbf{W}_V) + f_L(p_L^\ast,\mathbf{W}_L) + z\dbl v\dbr + \Delta u
    \end{matrix}\right).\label{2p_RP:ex_sys}
\end{align}
The quantities in the second component account for the classical outer waves, see Ref. \cite{Hantke2019a}. For $K\in\{V,L\}$ we have
    \begin{align*}
        f_K(p_K^\ast,\mathbf{W}_K) =
        \begin{dcases}
            \sqrt{-[\![p]\!][\![v]\!]},&\,\,p_K^\ast > p_K\,\,\text{(shock)},\\
            \int_{p_K}^{p_K^\ast} \frac{v_K(\zeta)}{a_K(\zeta)}\,\textup{d}\zeta,&\,\,p_K^\ast \leq p_K\,\,\text{(rarefaction)}.
        \end{dcases}
    \end{align*}
 The Jacobian is given by   
%
%\begin{widetext}
%
\begin{align}
    &\textbf{\textup{D}}\mathbf{G}(p_V^\ast,p_L^\ast) =\dots\notag\\
    &\left(\begin{matrix}
        -1 + 2z\dfrac{\partial z}{\partial p_V^\ast}\dbl v\dbr + z^2\left(\dfrac{v_V^\ast}{a_V^\ast}\right)^2 &
         1 + 2z\dfrac{\partial z}{\partial p_L^\ast}\dbl v\dbr - z^2\left(\dfrac{v_L^\ast}{a_L^\ast}\right)^2 \\[8pt]
        \dfrac{\dd}{\dd p_V^\ast}f_V(p_V^\ast,\mathbf{W}_V) + \dfrac{\partial z}{\partial p_V^\ast}\dbl v\dbr + z\left(\dfrac{v_V^\ast}{a_V^\ast}\right)^2 &
        \dfrac{\dd}{\dd p_L^\ast}f_L(p_L^\ast,\mathbf{W}_L) + \dfrac{\partial z}{\partial p_L^\ast}\dbl v\dbr - z\left(\dfrac{v_L^\ast}{a_L^\ast}\right)^2
    \end{matrix}\right).\label{jacobian:ex_sys}
\end{align}
%
%\end{widetext}

%The derivatives of the classical wave functions are given as in Subsection \ref{subsec:rp_sol_iso_gen} and the derivatives of the kinetic relation are as in Section \ref{sec:pb_sol_iso_gen}.\\
Once we have determined the pressures we can calculate the densities via the corresponding EOS, the mass flux $z$ via the kinetic relation,
the velocities via the wave relations for the classical waves, and the velocity of the phase boundary using the definition of the mass flux $z = -\rho(u - w)$.
Since we always align the computational grid with the phase boundary we have to use the corresponding flux when we apply \eqref{godunov_meth:nonunif}.
The flux at the phase boundary is determined using the jump conditions \eqref{jc:isoT:cons_mass_1d} and \eqref{jc:isoT:cons_mom_1d}. These may be rewritten as
\begin{align*}
    \rho_V^\ast(u_V^\ast - w) &= \rho_L^\ast(u_L^\ast - w) = -z,\\
    -z u_V^\ast + p_V^\ast &= -z u_L^\ast + p_L^\ast.
\end{align*}
Thus we directly read off the flux across the phase boundary, i.e.
\begin{align}
    \mathbf{F}_{PB} = \left[\begin{matrix} -z\\ -zu^* + p^*\end{matrix}\right].\label{pb_flux:pt}
\end{align}
Here one has the freedom to choose either the vapor or the liquid star state values.

So far we have discussed the case \emph{with} phase transition. The case \emph{without} phase transition is obtained using $z = 0$, which then implies $\dbl p\dbr = 0$ and $\dbl u\dbr = 0$.
Thus we have to solve only a single nonlinear equation, see Ref. \cite{Hantke2019a}. % given in Theorem \ref{thm:isoT:sol_without_pt}.
The flux is given by
\begin{equation}
    \mathbf{F}_{PB} = \left[\begin{matrix} 0\\ p^*\end{matrix}\right].\label{pb_flux:nopt}
\end{equation}
Further one immediately verifies that in this case the phase boundary is quite analogue to a contact discontinuity.
Therefore we use in our numerical simulations, for the case without phase transition, the wave speed estimates for the \emph{HLLC} solver as given in Refs. \cite{Batten1997,Toro2009}.
Given the (left) vapor state and the (right) liquid state we proceed as follows
\begin{align}
    S_L &= u_L + a_L,\notag\\
    S_V & = u_V - a_V,\label{HLLC:2p_nopt}\\
    w &= \frac{p_L - p_V + \rho_Vu_V(S_V - u_V) - \rho_Lu_L(S_L - u_L)}{\rho_V(S_V - u_V) - \rho_L(S_L - u_L)},\notag\\
    \rho_L^\ast &= \rho_L\frac{S_L - u_L}{S_L - w}.\notag
\end{align}
Here $S_V$ and $S_L$ denote the velocities of the classical waves. The density is calculated according to the HLLC solver as presented in Refs. \cite{Toro1994,Toro2009}.
The pressure $p^\ast$ may then be calculated using the liquid EOS. This procedure gives satisfactory results.

We use this simple calculation to choose proper initial pressure values for the Newton iteration in the case with phase transition.
It is shown in Ref. \cite{Hantke2019a} that the pressure in the star region for the case without phase transition correctly predicts the sign of the mass flux in the case with phase transition, i.e.\ whether we have evaporation or condensation.

\section{Phase Creation}\label{sec:num_cav_nucl}
The case of phase creation is a challenging issue: First, we have to detect that a new phase must be created.
%since we not only have to detect when phase creation occurs.
Then, we have to deal with technical problems such as the creation of new cells that are typically very small, potential interactions of multiple phase boundaries, and other problems.

%Further we also have to deal with technical problems like the creation of new cells which in most cases are very small and of course with multiple phase boundaries.
%
\subsection{Cavitation}
When we encounter a liquid/liquid Riemann problem we may observe cavitation, i.e.,\ the creation of vapor. For the detailed analysis we again refer to Refs. \cite{Thein2018,Hantke2019a}. 
%is presented in Subsection \ref{subsec:cav_iso_gen}.
In view of the results given there %in Theorem \ref{thm:isoT:sol_rp_1p_liq}, Definition \ref{def:isoT:cav_crit} and Theorem \ref{thm:isoT:evap_exp}
the main outline for the numerics is as follows.
We solve the single phase Riemann problem in the liquid phase.
If there is no solution to this problem, in particular when the liquid pressure in the star region is smaller than the predefined minimum liquid pressure, $p^\ast < p_{min}$, we have cavitation.
In this case we store the position of the involved cells and perform an extra calculation. We summarize the performed steps:
\begin{enumerate}[(i)]
    \item Solve the single phase Riemann problem without phase transition; if $p^\ast < p_{min}$ cavitation occurs.
    \item Solve the single phase Riemann problem with phase transition.% according to Theorem \ref{thm:isoT:evap_exp}.%\footnote{For the solution one may also use the approximation presented
    %in Section \ref{sec:approx_kinrel}}
    \item From the solution we obtain $p_V^\ast,\,z,\,w_{left}$ and $w_{right}$.
    \item Create a vapor cell of size $(w_{right} - w_{left})\Delta t$ with the cell values $\rho_V^\ast,\,p_V^\ast$ and $u_V^\ast$.
    \item The fluxes at the phase boundaries are given by
    \begin{align}
        \mathbf{F}_{PB}^{(left)} = \left[\begin{matrix} z\\ -zu_V^* + p_V^*\end{matrix}\right]
        \quad\text{and}\quad
        \mathbf{F}_{PB}^{(right)} = \left[\begin{matrix} -z\\ -zu_V^* + p_V^*\end{matrix}\right].\label{flux:cavitation}
    \end{align}
\end{enumerate}
In the next time step we then have two phase boundaries which are treated as discussed above.%in Section \ref{sec:num_interface}.
\subsection{Nucleation}
The case of nucleation is treated analogous to the previous case of cavitation.
When we encounter a vapor/vapor Riemann problem we may observe nucleation, i.e.\ the creation of liquid. The detailed analysis is presented in Refs. \cite{Thein2018,Hantke2019a}. %Subsection \ref{subsec:nucl_iso_gen}.
In view of the results given there %in Theorem \ref{thm:isoT:sol_rp_1p_vap}, Definition \ref{def:isoT:nucl_crit} and Theorem \ref{thm:isoT:cond_compr}
the main outline for the numerics is as follows.
We solve the single phase Riemann problem in the vapor phase.
If there is no solution to this problem, i.e.\ when the vapor pressure in the star region is greater than the predefined maximum vapor pressure, $p^\ast > \tilde{p}$, we have nucleation.
In this case we store the position of the involved cells and perform an extra calculation. We summarize the performed steps:
\begin{enumerate}[(i)]
    \item Solve the single phase Riemann problem without phase transition; if $p^\ast > \tilde{p}$ nucleation occurs.
    \item Solve the single phase Riemann problem with phase transition.% according to Theorem \ref{thm:isoT:cond_compr}
    %\footnote{For the solution one may also use the approximation presented in Section\ref{sec:approx_kinrel}}.
    \item From the solution we obtain $p_L^\ast,\,z,\,w_{left}$ and $w_{right}$.
    \item Create a liquid cell of size $(w_{right} - w_{left})\Delta t$ with the cell values $\rho_L^\ast,\,p_L^\ast$ and $u_L^\ast$.
    \item The fluxes at the phase boundaries are given by
    \begin{align}
        \mathbf{F}_{PB}^{(left)} = \left[\begin{matrix} -z\\ -zu_L^* + p_L^*\end{matrix}\right]
        \quad\text{and}\quad
        \mathbf{F}_{PB}^{(right)} = \left[\begin{matrix} z\\ -zu_L^* + p_L^*\end{matrix}\right].\label{flux:nucleation}
    \end{align}
\end{enumerate}
Again, in the next time step we then have two phase boundaries which are treated as discussed above.%in Section \ref{sec:num_interface}.
%Especially in the case of nucleation it is important to not that the created cells are usually very small (e.g.\ seven orders of magnitude) compared %to the other cells.
%Thus the local time stepping presented in Section \ref{sec:num_interface} is very much needed to avoid too small time steps in the whole %computational domain.
%We will again comment on this in the following Section \ref{sec:num_ex}